\documentclass{svmult}

\usepackage{amsmath}
\usepackage{amssymb}
\usepackage{epsfig}

\def\CC{{\mathbb C}}

\def\HH{{\mathbb H}}

\def\NN{{\mathbb N}}

\def\RR{{\mathbb R}}
\def\SSS{{\mathbb S}}

\def\ZZ{{\mathbb Z}}

\def\GamHH{\Gamma\backslash{\mathbb H}^2}
\def\ZHH{\Z\backslash{\mathbb H}^2}
\def\FGam{{\mathcal F}_{\Gamma}}

\def\D{{\mathfrak D}}
\def\H{{\mathfrak H}}
\def\Z{{\mathfrak Z}}
\def\Isom{\operatorname{Isom}}
\def\e{\mathrm{e}}
\def\i{\mathrm{i}}

\def\C{\operatorname{C{}}}

\def\L{\operatorname{L{}}}

\def\res{\operatorname{res}}

\def\SL{\operatorname{SL}}

\def\PSL{\operatorname{PSL}}

\def\U{\operatorname{U{}}}

\def\tr{\operatorname{tr}}
\def\Tr{\operatorname{Tr}}

\def\supp{\operatorname{supp}}
\def\Vol{\operatorname{Vol}}
\def\Area{\operatorname{Area}}

\def\vecm{{\text{\boldmath$m$}}}
\def\vecn{{\text{\boldmath$n$}}}

\def\vecx{{\text{\boldmath$x$}}}

\def\vecrho{{\text{\boldmath$\rho$}}}

\def\vectau{{\text{\boldmath$\tau$}}}

\def\scrC{{\mathcal C}}
\def\scrD{{\mathcal D}}
\def\scrF{{\mathcal F}}

\def\scrM{{\mathcal M}}

\def\Re{\operatorname{Re}}
\def\Im{\operatorname{Im}}

\begin{document}

\title*{Selberg's trace formula: an introduction}

\author{Jens Marklof}
\institute{School of Mathematics, University of Bristol,
Bristol BS8 1TW, U.K.
\texttt{j.marklof@bristol.ac.uk}}

\maketitle

The aim of this short lecture course is to develop Selberg's trace formula 
for a compact hyperbolic surface $\scrM$, and discuss some of its 
applications. The main motivation for our studies is {\em quantum chaos}:
the Laplace-Beltrami operator $-\Delta$ on the surface 
$\scrM$ represents the quantum Hamiltonian of a particle,
whose classical dynamics is governed by the (strongly chaotic) geodesic flow
on the unit tangent bundle of $\scrM$. 
The trace formula is currently the only available tool to analyze
the fine structure of the spectrum of $-\Delta$; no individual formulas for 
its eigenvalues are known. In the case of more general quantum systems, 
the role of Selberg's formula is taken over by the semiclassical
{\em Gutzwiller trace formula} \cite{Gutzwiller89}, \cite{Combescure99}.

We begin by reviewing the trace formulas for the simplest compact manifolds, 
the circle $\SSS^1$ (Section \ref{secPoisson}) and the sphere $\SSS^2$
(Section \ref{secTrace}). In both cases, the corresponding geodesic
flow is integrable, and the trace formula
is a consequence of the {\em Poisson summation formula}.
In the remaining sections we shall discuss the following topics:
{\em the Laplacian on the hyperbolic plane and isometries}
(Section \ref{secHyperbolic}); 
{\em Green's functions} (Section \ref{secGreen});
{\em Selberg's point pair invariants} (Section \ref{secSelberg});
{\em The ghost of the sphere} (Section \ref{ghost});
{\em Linear operators on hyperbolic surfaces} (Section \ref{secHyperbolic2});
{\em A trace formula for hyperbolic cylinders and
poles of the scattering matrix} (Section \ref{secCylinders});
{\em Back to general hyperbolic surfaces} (Section \ref{secBack});
{\em The spectrum of a compact surface, 
Selberg's pre-trace and trace formulas} (Section \ref{secSpectrum});
{\em Heat kernel and Weyl's law} (Section \ref{secWeyl});
{\em Density of closed geodesics} (Section \ref{secHuber});
{\em Trace of the resolvent} (Section \ref{secResolvent});
{\em Selberg's zeta function} (Section \ref{selzeta});
{\em Suggestions for exercises and further reading} 
(Section \ref{secSuggestions}).

Our main references are Hejhal's classic lecture notes 
\cite[Chapters {\sc one} and {\sc two}]{HejhalI}, Balazs and Voros'
excellent introduction \cite{Balazs86},
and Cartier and Voros' {\em nouvelle interpr\'etation} \cite{Cartier90}.
Section \ref{secSuggestions} comprises a list of references 
for further reading.

These notes are based on lectures given at the 
International School {\em Quantum Chaos on Hyperbolic Manifolds},
Schloss Reisensburg (G\"unzburg, Germany), 4-11 October 2003.

\section{Poisson summation \label{secPoisson}}

The Poisson summation formula reads
\begin{equation}\label{Poisson}
\sum_{m\in\ZZ} h(m)= \sum_{n\in\ZZ} \int_\RR h(\rho) \, \e^{2\pi\i n \rho}\,
d\rho
\end{equation}
for any sufficiently nice test function $h:\RR\to\CC$.
One may for instance take $h\in\C^2(\RR)$ with 
$|h(\rho)|\ll (1+|\rho|)^{-1-\delta}$ for some $\delta>0$. 
(The notation $x\ll y$ means here {\em there exists a constant $C>0$ such that
$x \leq C y$}.)
Then
both sums in \eqref{Poisson} converge absolutely.
\eqref{Poisson} is proved by expanding the periodic
function
\begin{equation}
f(\rho)=\sum_{m\in\ZZ} h(\rho+m)
\end{equation}
in its Fourier series, and then setting $\rho=0$.

The Poisson summation formula is our first example of a {\em trace formula}:
The eigenvalues of the positive Laplacian $-\Delta=-\frac{d^2}{dx^2}$
on the circle $\SSS^1$ of length $2\pi$ are $m^2$ where $m=0,\pm1,\pm2,\ldots$,
with corresponding eigenfunctions $\varphi_m(x)=(2\pi)^{-1/2}\e^{\i m x}$.
Consider the linear operator $L$ acting on $2\pi$-periodic functions by
\begin{equation}
[L f](x):= \int_0^{2\pi} k(x,y) f(y)\, dy
\end{equation}
with kernel
\begin{equation}
k(x,y)= \sum_{m\in\ZZ} h(m)\, \varphi_m(x)\, \overline{\varphi}_m(y) .
\end{equation}
Then 
\begin{equation}
L \varphi_m = h(m) \varphi_m
\end{equation}
and hence the Poisson summation formula says that
\begin{equation}
\Tr L = \sum_{m\in\ZZ} h(m)= 
\sum_{n\in\ZZ} \int_\RR h(\rho) \, \e^{2\pi\i n \rho}\,
d\rho .
\end{equation}
The right hand side in turn has a geometric interpretation as a sum
over the periodic orbits of the geodesic flow on $\SSS^1$, whose
lengths are $2\pi |n|$, $n\in\ZZ$.

An important example for a linear operator of the above type 
is the resolvent of the Laplacian,
$(\Delta+\rho^2)^{-1}$, with $\Im\rho<0$.
The corresponding test function is 
$h(\rho')= (\rho^2-{\rho'}^2)^{-1}$. Poisson summation yields in this case
\begin{equation}\label{Poisson2}
\sum_{m\in\ZZ} (\rho^2-m^2)^{-1}
= \sum_{n\in\ZZ} \int_\RR \frac{\e^{-2\pi\i |n| \rho'}}{\rho^2-{\rho'}^2}\,
d\rho'
\end{equation}
and by shifting the contour to $-\i\infty$ and collecting the residue
at $\rho'=\rho$,
\begin{center}
\epsfig{figure=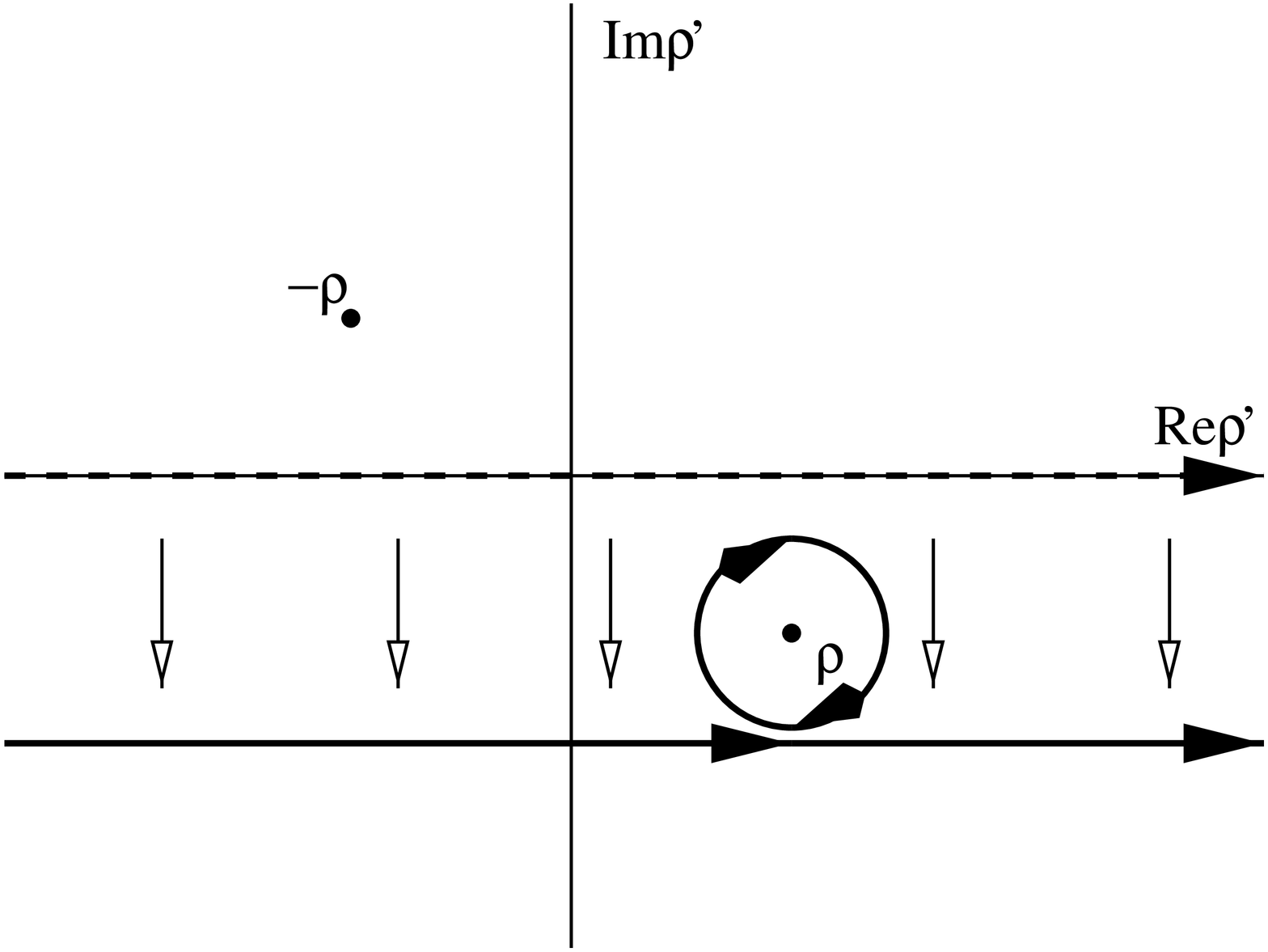,width=0.35\textwidth,angle=270}
\end{center}
we find
\begin{equation}\label{Poisson3}
\sum_{m\in\ZZ} (\rho^2-m^2)^{-1}
= \frac{\pi\i}{\rho} \sum_{n\in\ZZ} \e^{-2\pi\i |n| \rho} .
\end{equation}
The right hand side resembles
the geometric series expansion of $\cot z$ for $\Im z<0$,
\begin{equation}\label{cot}
\cot z = \frac{2\i \e^{-\i z}\cos z}{1-e^{-2\i z}}
= \i (1+\e^{-2\i z}) \sum_{n=0}^\infty \e^{-2\i n z}
= \i \sum_{n\in\ZZ} \e^{-2\i |n| z} .
\end{equation}
Hence
\begin{equation}\label{Poisson4}
\sum_{m\in\ZZ} (\rho^2-m^2)^{-1}
= \frac{\pi}{\rho} \cot(\pi\rho),
\end{equation}
which can also be written in the form
\begin{equation}\label{although}
\frac12 \sum_{m\in\ZZ} \left[\frac{1}{\rho-m}+\frac{1}{\rho+m} \right]
= \pi \cot(\pi\rho),
\end{equation}
The above $h$ is an example of a test function with particularly
useful analytical properties. More generally, suppose 
\begin{enumerate}
\item[(i)]
$h$ is analytic for $|\Im\rho|\leq \sigma$ for some $\sigma>0$;
\item[(ii)]
$|h(\rho)|\ll(1+|\Re\rho|)^{-1-\delta}$ for some $\delta>0$,
uniformly for all $\rho$ in the strip $|\Im\rho|\leq \sigma$.
\end{enumerate}

\begin{theorem}\label{circthm}
If $h$ satisfies {\rm (i), (ii)}, then
\begin{equation} \label{circeq}
\sum_{m\in\ZZ} h(m)
=
\frac{1}{2\i} \int_{\scrC_=}  h(\rho) \, \cot(\pi\rho)\, d\rho 
\end{equation}
where the path of integration $\scrC_=$ is
\begin{center}
\epsfig{figure=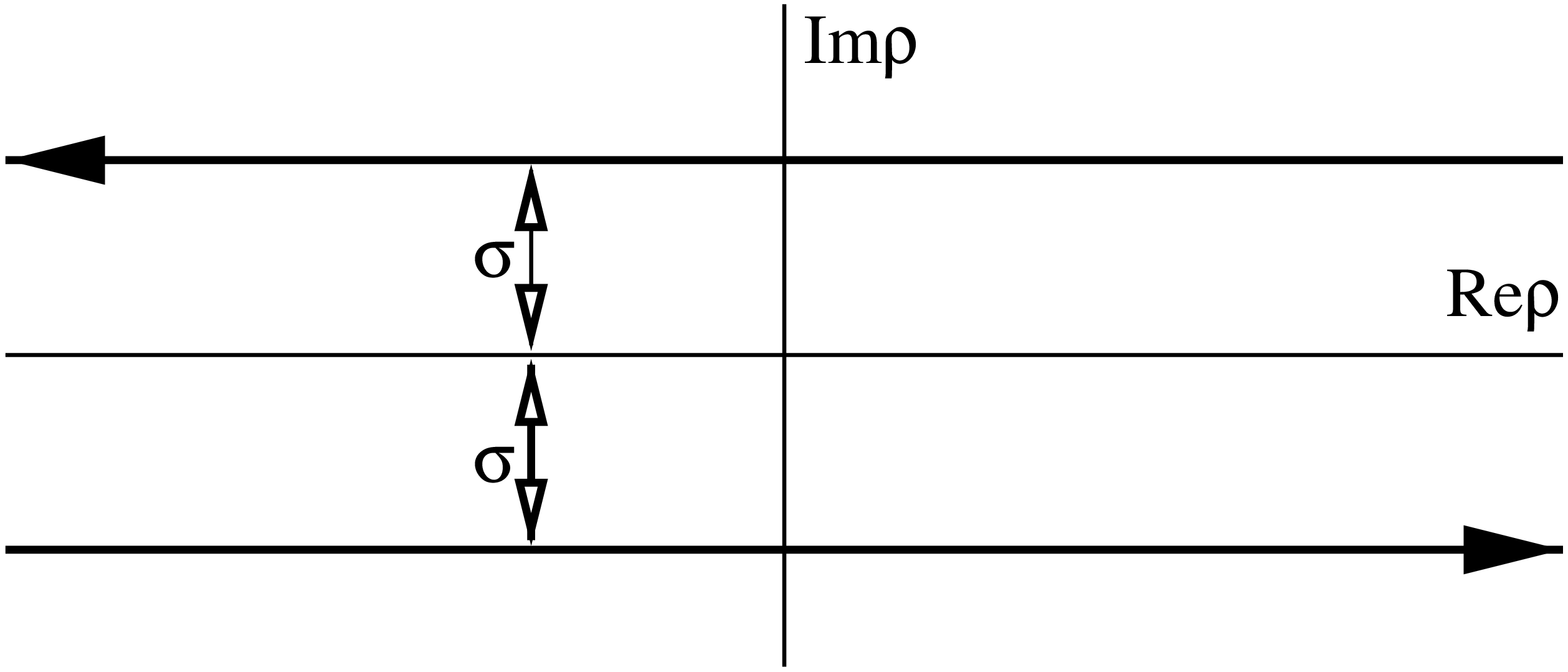,width=0.2\textwidth,angle=270}
\end{center}
\end{theorem}

\begin{proof}
The Poisson summation formula \eqref{Poisson} may be written in the form
\begin{equation}\label{Poisson20}
\sum_{m\in\ZZ} h(m) 
=\frac12 \sum_{n\in\ZZ} \int_\RR [h(\rho)+h(-\rho)] \, 
\e^{-2\pi\i |n| \rho}\,d\rho .
\end{equation}
We shift the contour of the integral to 
$\int_{-\infty-\i\sigma}^{\infty-\i\sigma}$.
The geometric series expansion of $\cot z$ in \eqref{cot}
converges absolutely, uniformly for all $z$ with fixed negative imaginary part.
We may therefore exchange summation and integration, 
\begin{equation}
\sum_{m\in\ZZ} h(m) =
\frac{1}{2\i} \int_{-\infty-\i\sigma}^{\infty-\i\sigma}  
[h(\rho)+h(-\rho)] \, \cot(\pi\rho)\, d\rho .
\end{equation}
We conclude with the observation that
\begin{equation}
\frac{1}{2\i} \int_{-\infty-\i\sigma}^{\infty-\i\sigma}  
h(-\rho) \, \cot(\pi\rho)\, d\rho
=
\frac{1}{2\i} \int_{\infty+\i\sigma}^{-\infty+\i\sigma}  
h(\rho) \, \cot(\pi\rho)\, d\rho
\end{equation}
since $\cot z$ is odd.
\qed
\end{proof}

\begin{remark}
This theorem can of course also be proved by shifting the lower 
contour in \eqref{circeq} across the poles to the upper contour, 
and evaluating the corresponding residues.
\end{remark}

\section{A trace formula for the sphere \label{secTrace}}

The Laplacian on the sphere $\SSS^2$ is given by
\begin{equation}
\Delta= 
\frac{1}{\sin\theta} \frac{\partial}{\partial \theta}
\left(\sin\theta  \frac{\partial}{\partial \theta}\right) 
+ \frac{1}{\sin^2\theta} \frac{\partial^2}{\partial \phi^2}
\end{equation}
where $\theta\in[0,\pi)$, $\phi\in[0,2\pi)$ are the standard 
spherical coordinates. The eigenvalue problem
\begin{equation}
(\Delta+\lambda) f = 0
\end{equation}
is solved by the spherical harmonics $f=Y_l^m$
for integers $l=0,1,2,\ldots$, $m=0,\pm1,\pm2,\ldots,\pm l$,
where
\begin{equation}
Y_l^m(\theta,\phi)= (-1)^m \left[\frac{(2l+1)}{4\pi}\,
\frac{(l-m)!}{(l+m)!}\right]^{1/2} P_l^m(\cos\theta)\, \e^{\i m\phi}
\end{equation}
and $P_l^m$ denotes the associated Legendre function of the first kind.
The eigenvalue corresponding to $Y_l^m$ is $\lambda = l(l+1)$,
and hence appears with multiplicity $2l+1$. Let us label
all eigenvalues (counting multiplicity) by
\begin{equation}
0=\lambda_0 < \lambda_1 \leq \lambda_2 \leq \ldots \to\infty ,
\end{equation}
and set $\rho_j=\sqrt{\lambda_j+\tfrac 14}>0$.
For any even test function $h\in\C^2(\RR)$ with the bound 
$|h(\rho)|\ll(1+|\Re\rho|)^{-2-\delta}$ for some $\delta>0$
(assume this bound also holds for the first and second derivative)
we have
\begin{align}
\sum_{j=0}^\infty h(\rho_j)
& =
\sum_{l=0}^\infty (2l+1)\, h(l+\tfrac12) \\
& =
\sum_{l=-\infty}^\infty |l+\tfrac12|\, h(l+\tfrac12)\\
& =
\sum_{n=-\infty}^\infty \int_\RR |l+\tfrac12|\, h(l+\tfrac12)
\e^{2\pi\i l n} dl \\
& = \label{sun1}
2 \sum_{n\in\ZZ} 
(-1)^n
\int_0^\infty \rho \,h(\rho)\, \cos(2\pi n \rho)\, d\rho ,
\end{align}
in view of the Poisson summation formula.
We used the test function $|\rho| h(\rho)$ which is not 
continuously differentiable at $\rho=0$. This is not a problem, since
we check that (using integration by parts twice) 
\begin{equation}
\int_0^\infty \rho \,h(\rho)\, \cos(2\pi n \rho)\, d\rho=O(n^{-2})
\end{equation}
hence all sums are absolutely convergent. 
With $\Area(\SSS^2)=4\pi$,
it is suggestive to write the trace formula \eqref{sun1} for the sphere 
in the form
\begin{equation}
\sum_{j=0}^\infty h(\rho_j)
=
\frac{\Area(\SSS^2)}{4\pi}
\int_\RR |\rho| \,h(\rho)\, d\rho 
+\sum_{n\neq 0} 
(-1)^n
\int_\RR |\rho| \,h(\rho)\, \e^{2\pi\i n \rho}\, d\rho .
\end{equation}
As in the trace formula for the circle, the sum on the right
hand side may again be viewed as a sum over the closed geodesics
of the sphere which, of course, all have lengths $2\pi|n|$.
The factor $(-1)^n$ accounts for the number of conjugate points
traversed by the corresponding orbit.

The sum in \eqref{sun1} resembles the geometric series
expansion for $\tan z$ for $\Im z<0$,
\begin{equation}\label{tan}
\tan z = - \cot(z+\pi/2)
= - \i \sum_{n\in\ZZ} (-1)^n \e^{-2\i |n| z} .
\end{equation}
As remarked earlier, the sum converges uniformly for all $z$ with
fixed $\Im z < 0$. We have in fact the uniform bound
\begin{equation}\label{unitan}
\sum_{n\in\ZZ} \left| (-1)^n \e^{-2\i |n| z} \right|
\leq
1+2 \sum_{n=1}^\infty \e^{2 n\Im z}
\leq
1+ 2 \int_0^\infty \e^{2 x \Im z} \, dx
= 1 - \frac{1}{\Im z}
\end{equation}
which holds for all $z$ with $\Im z<0$.

Let us use this identity to rewrite the trace formula.
Assume $h$ satisfies the following hypotheses.
\begin{enumerate}
\item[(i)]
$h$ is analytic for $|\Im\rho|\leq \sigma$ for some $\sigma>0$;
\item[(ii)] 
$h$ is even, i.e., $h(-\rho)=h(\rho)$;
\item[(iii)]
$|h(\rho)|\ll(1+|\Re\rho|)^{-2-\delta}$ for some $\delta>0$,
uniformly for all $\rho$ in the strip $|\Im\rho|\leq \sigma$.
\end{enumerate}

\begin{theorem}
If $h$ satisfies {\rm (i), (ii), (iii)}, then
\begin{equation} \label{sph}
\sum_{j=0}^\infty h(\rho_j)
= - \frac{1}{2\i} \int_{\scrC_\times} h(\rho)\,\rho \tan(\pi\rho)\,d\rho ,
\end{equation}
with the path of integration
\begin{center}
\epsfig{figure=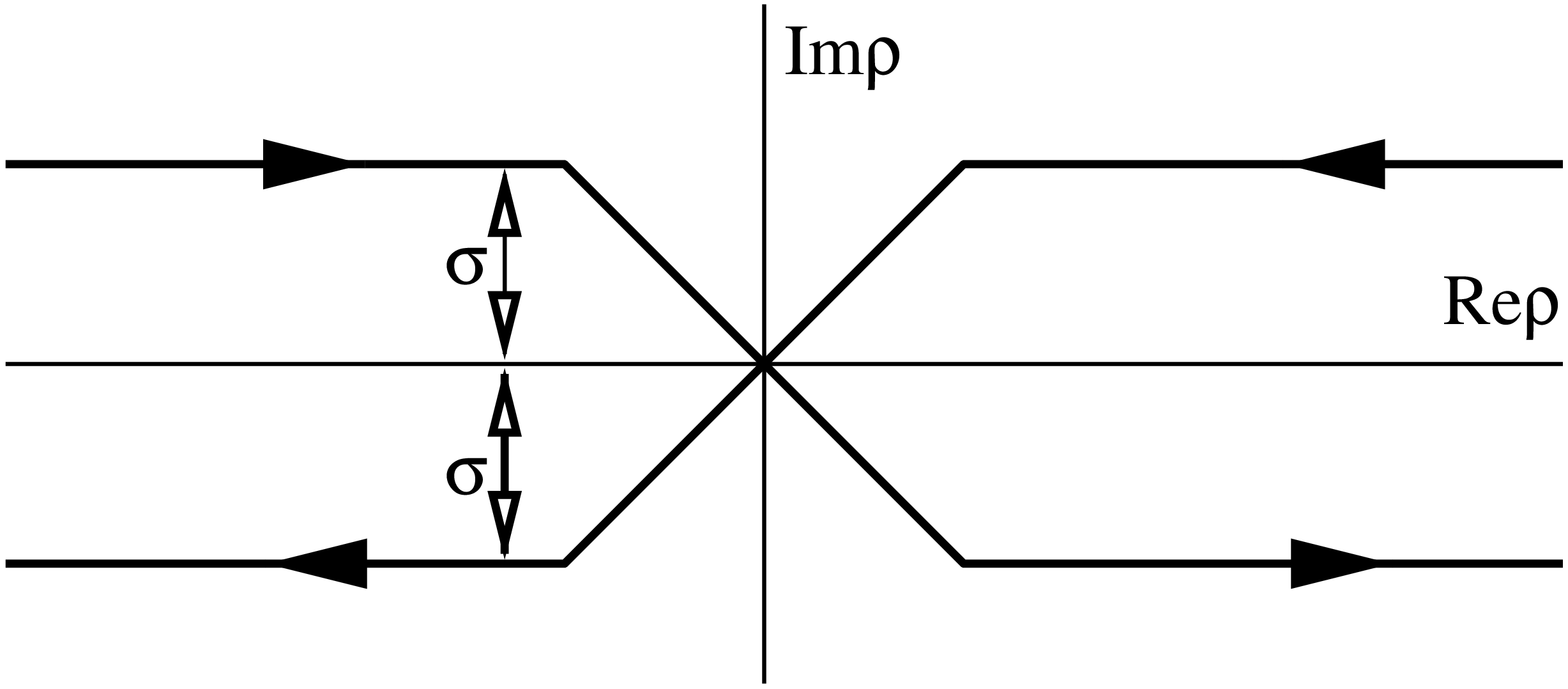,width=0.2\textwidth,angle=270}
\end{center}
\end{theorem}

\begin{proof}
We express \eqref{sun1} in the form
\begin{equation}
\sum_{n\in\ZZ} 
(-1)^n
\int_0^\infty
\rho \,h(\rho)\, \e^{2\pi\i |n| \rho}\, d\rho 
+
\sum_{n\in\ZZ} 
(-1)^n
\int_0^\infty
\rho \,h(\rho)\, \e^{-2\pi\i |n| \rho}\, d\rho .
\end{equation}
which equals
\begin{equation}\label{onetwo}
- \sum_{n\in\ZZ} 
(-1)^n
\int_{-\infty}^0
\rho \,h(\rho)\, \e^{-2\pi\i |n| \rho}\, d\rho 
+ \sum_{n\in\ZZ} 
(-1)^n
\int_0^\infty
\rho \,h(\rho)\, \e^{-2\pi\i |n| \rho}\, d\rho .
\end{equation}
Let us first consider the second integral.
We change the path of integration to $\scrC_1$:
\begin{center}
\epsfig{figure=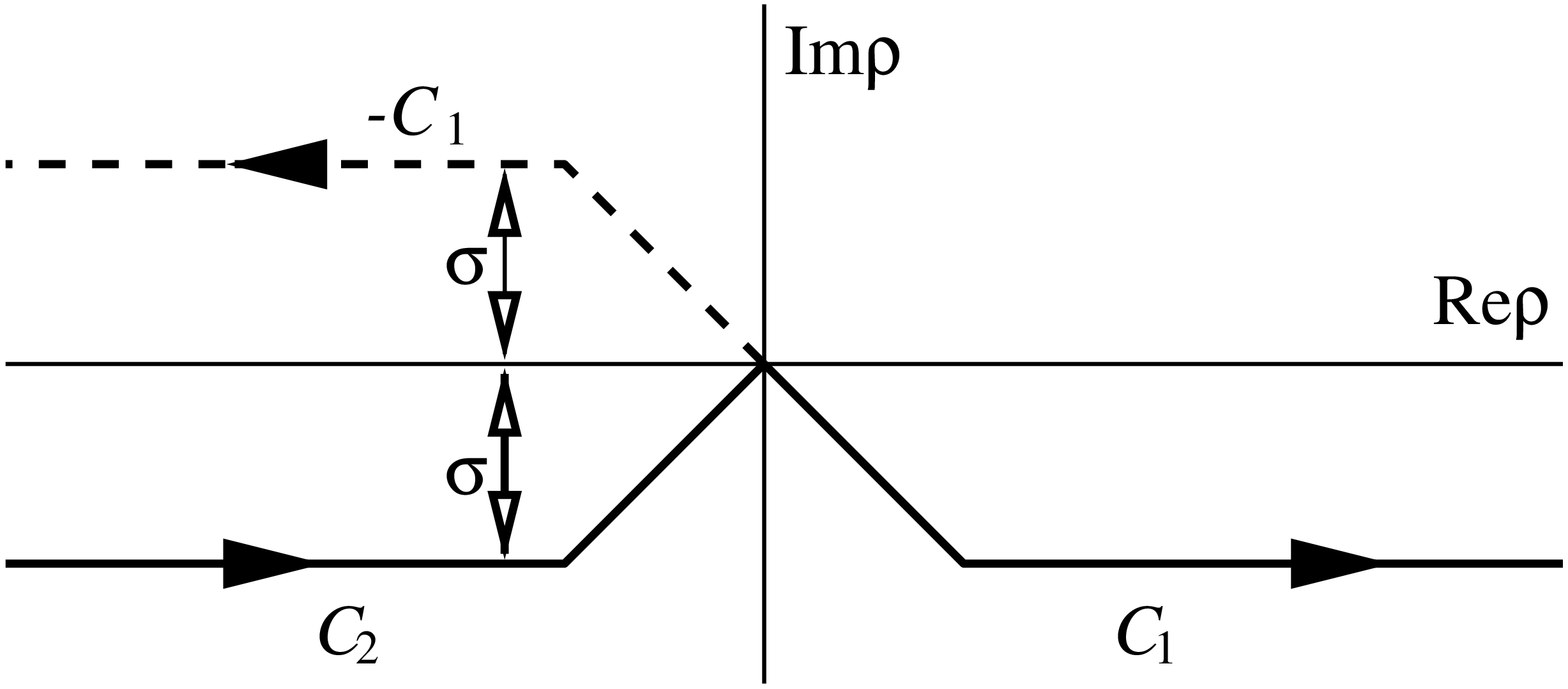,width=0.2\textwidth,angle=270}
\end{center}
Due to the uniform bound
\eqref{unitan} and
\begin{equation}
\int_{\scrC_1}
\left| \rho \,h(\rho)[1 - (2\pi\Im \rho)^{-1}]\,  d\rho \right| < \infty
\end{equation} 
we may exchange integration and summation, and hence the second integral
in \eqref{onetwo} evaluates to
\begin{equation}
\sum_{n\in\ZZ} 
(-1)^n
\int_0^{\infty}
\rho \,h(\rho)\, \e^{-2\pi\i |n| \rho}\, d\rho 
=
\i \int_{\scrC_1}
h(\rho)\,\rho \tan(\pi\rho)\,d\rho 
\end{equation}
The first integral in \eqref{onetwo} is analogous, we have
\begin{align}
- \sum_{n\in\ZZ} 
(-1)^n
\int_{-\infty}^0
\rho \,h(\rho)\, \e^{-2\pi\i |n| \rho}\, d\rho 
& =
-\i \int_{\scrC_2}
h(\rho)\,\rho \tan(\pi\rho)\,d\rho \\
& =
\i \int_{\scrC_2^{-1}}
h(\rho)\,\rho \tan(\pi\rho)\,d\rho .
\end{align}
The final result is obtained by reflecting these paths at the origin,
using the fact that $h$ is even.
\qed
\end{proof}

\begin{remark}
The poles of $\tan z$
and corresponding residues can be easily worked out from \eqref{although},
\begin{align}
\pi \tan(\pi\rho) & = -\pi \cot\left[\pi(\rho+\tfrac12)\right] \\
& = 
- \frac{1}{2}\sum_{l=-\infty}^\infty  
\left[ \frac{1}{\rho-(l-\tfrac12)}+ \frac{1}{\rho+(l+\tfrac12)} \right] 
\label{here12} \\
& = 
-\frac{1}{2}\sum_{l=-\infty}^\infty  
\left[\frac{1}{\rho+(l+\tfrac12)} + \frac{1}{\rho-(l+\tfrac12)} \right] 
\label{here13} \\
\intertext{(the sum has not been reordered, we
have simply shifted the bracket)}
& = 
- \sum_{l=0}^\infty  
\left[\frac{1}{\rho+(l+\tfrac12)} + \frac{1}{\rho-(l+\tfrac12)} \right] .
\label{here14} 
\end{align}
Note that the extra factor $\rho$ in the integral \eqref{sph}, as compared
to Theorem \ref{circthm}, yields the multiplicity of the eigenvalues of the
sphere.
\end{remark}

\section{The hyperbolic plane \label{secHyperbolic}}

In this section we briefly summarize some basic features
of hyperbolic geometry; for a detailed discussion see 
Buser's lecture notes \cite{Buser}.

The hyperbolic plane $\HH^2$ may be abstractly defined as 
the simply connected two-dimensional Riemannian manifold
with Gaussian curvature $-1$. Let us introduce three convenient
coordinate systems for $\HH^2$: 
the Poincare disk $\D=\{ z : |z|<1 \}$, the upper half plane
$\H=\{ z : \Im z >0\}$ 
and polar coordinates $(\tau,\phi)\in\RR_{\geq 0}\times [0,2\pi)$.
In these parametrizations,
the line element $ds$, volume element $d\mu$ and the Riemannian distance 
$d(z,z')$ between two points $z,z'\in\HH^2$ are as follows.
\begin{center}
\begin{tabular}{|c|c|c|c|}
\hline
$\HH^2$ & $ds^2$ & $d\mu$ & $\cosh d(z,z')$ \\ 
\hline
$\D$  & $\displaystyle \frac{4(dx^2+dy^2)}{(1-x^2-y^2)^2}$ & 
$\displaystyle \frac{4\,dx\, dy}{(1-x^2-y^2)^2}$ &
$\displaystyle 1 + \frac{2|z-z'|^2}{(1-|z|)^2(1-|z'|)^2}$ \\
\hline
$\H$ & $\displaystyle \frac{dx^2+dy^2}{y^2}$ & 
$\displaystyle \frac{dx\,dy}{y^2}$ & 
$\displaystyle 1 + \frac{|z-z'|^2}{2\Im z \Im z'}$  \\ 
\hline
polar & $d\tau^2+\sinh^2\tau \, d\phi^2$ & 
$\sinh\tau\, d\tau\, d\phi$ & 
$\cosh\tau$ [for $z=(\tau,\phi)$, $z'=(0,0)$]\\ 
\hline 
\end{tabular}
\end{center}
The {\em group of isometries} of $\HH^2$, denoted by $\Isom(\HH^2)$,
is the group of smooth coordinate transformations which leave the
Riemannian metric invariant. The group of {\em orientation preserving}
isometries is called $\Isom^+(\HH^2)$.
We define the {\em length} of an isometry
$g\in\Isom(\HH^2)$ by
\begin{equation}
\ell_g=\ell(g) = \inf_{z\in\HH^2} d(gz,z) .
\end{equation}
Those $g\in\Isom^+(\HH^2)$ for which $\ell>0$ are called {\em hyperbolic}.
In the half plane model, $\Isom^+(\HH^2)$ acts by fractional linear
transformations,
\begin{equation}
g: \H \to \H, \qquad z \mapsto gz:= \frac{az+b}{cz+d} , \qquad 
\begin{pmatrix} a & b \\ c & d \end{pmatrix}\in\SL(2,\RR) 
\end{equation}
(we only consider orientation-preserving isometries here).
We may therefore identify $g$ with a matrix in $\SL(2,\RR)$,
where the matrices $g$ and $-g$ obviously correspond to the 
same fractional linear transformation. $\Isom^+(\HH^2)$ may thus
be identified with the group $\PSL(2,\RR)=\SL(2,\RR)/\{\pm 1\}$.
In this representation,
\begin{equation}
2\cosh\big(\ell_g/2\big) = \max\{|\tr g|, 2\} ,
\end{equation}
since every matrix $g\in\SL(2,\RR)$ is conjugate to one of the following
three,
\begin{equation}
\begin{pmatrix} 1 & b \\ 0 & 1 \end{pmatrix}, \qquad
\begin{pmatrix} a & 0 \\ 0 & a^{-1} \end{pmatrix}, \qquad
\begin{pmatrix} \cos(\theta/2) & \sin(\theta/2) \\ 
-\sin(\theta/2) & \cos(\theta/2) \end{pmatrix} ,
\end{equation}
with $b\in\RR$, $a\in\RR_{>0}$, $\theta\in[0,2\pi)$.

The {\em Laplace-Beltrami operator} (or {\em Laplacian} for short) $\Delta$ 
of a smooth Riemannian manifold with metric
\begin{equation}
ds^2 = \sum_{ij} g_{jk} dx^j dx^k
\end{equation}
is given by the formula
\begin{equation}
\Delta= \sum_{ij} \frac{1}{\sqrt{g}} \frac{\partial}{\partial x^j}
\bigg( \sqrt{g}\, g^{jk} \frac{\partial}{\partial x^k} \bigg)
\end{equation}
where $g^{jk}$ are the matrix coefficients of the the inverse of the
matrix $(g_{jk})$, and $g=|\det(g_{jk})|$. In the above coordinate systems
for $\HH^2$ the Laplacian takes the following form.
\begin{center}
\begin{tabular}{|c|c|c|c|}
\hline
$\HH^2$ & $\Delta$ \\ 
\hline
$\D$  & 
$\displaystyle \frac{(1-x^2-y^2)^2}{4} 
\left( \frac{\partial^2}{\partial x^2} + \frac{\partial^2}{\partial y^2} 
\right)$ \\
\hline
$\H$ & 
$\displaystyle y^2 
\left( \frac{\partial^2}{\partial x^2} + \frac{\partial^2}{\partial y^2} 
\right)$ \\ 
\hline
polar &
$\displaystyle\frac{1}{\sinh\tau} \frac{\partial}{\partial \tau}
\left(\sinh\tau  \frac{\partial}{\partial \tau}\right) 
+ \frac{1}{\sinh^2\tau} \frac{\partial^2}{\partial \phi^2}$ \\ 
\hline 
\end{tabular}
\end{center}
One of the important properties of the Laplacian is that 
it commutes with every isometry $g\in\Isom(\HH^2)$. That is, 
\begin{equation} \label{commu}
\Delta T_g = T_g \Delta \qquad \forall g\in\Isom(\HH^2) .
\end{equation}
where the left translation operator $T_g$ acting on functions $f$ on $\HH^2$
is defined by
\begin{equation}\label{trans}
[T_g f](z) = f(g^{-1}z) 
\end{equation}
with $g\in\Isom(\HH^2)$. Even though \eqref{commu} is an intrinsic property
of the Laplacian and is directly related to the invariance of the Riemannian
metric under isometries, it is a useful exercise to verify \eqref{commu}
explicitly. To this end note that every isometry
may be represented as a product of fractional linear transformations
of the form
$z \mapsto a z$ $(a>0)$, $z \mapsto z+b$ $(b\in\RR)$,
$z \mapsto -1/z$, $z \mapsto -\overline z$.
It is therefore sufficient to check \eqref{commu} 
only for these four transformations.

\section{Green's functions \label{secGreen}}

The Green's function $G(z,w;\lambda)$ 
corresponding to the differential equation 
$(\Delta+\lambda) f(z) =0$ is formally defined as the integral kernel
of the resolvent $(\Delta+\lambda)^{-1}$, i.e., by the equation
\begin{equation}\label{greendiff0}
(\Delta+\lambda)^{-1} f(z) = \int G(z,w;\lambda)\, f(w)\, d\mu(w)
\end{equation}
for a suitable class of test functions $f$. A more precise characterization
is as follows:
\begin{description}
\item[(G1)]
$G(\,\cdot\,,w;\lambda)\in \C^\infty(\HH^2-\{w\})$ for every fixed $w$; 
\item[(G2)] 
$(\Delta+\lambda) G(z,w;\lambda) = \delta(z,w)$
for every fixed $w$; 
\item[(G3)]
as a function of $(z,w)$, $G(z,w;\lambda)$ depends on
the distance $d(z,w)$ only;
\item[(G4)]
$G(z,w;\lambda)\to 0$ as $d(z,w)\to\infty$. 
\end{description}
Here $\delta(z,w)$ is the Dirac distribution at $w$ with respect to the
measure $d\mu$. It is defined by the properties that 
\begin{description}
\item[(D1)]
$\delta(z,w) d\mu(z)$ is a probability measure on $\HH^2$;
\item[(D2)]
$\int_{\HH^2} f(z)\,
\delta(z,w)\, d\mu(z) = f(w)$ for all $f\in\C(\HH^2)$.
\end{description}
E.g., in the disk coordinates $z=x+\i y$,
$w=u+\i v\in\D$ we then have 
\begin{equation}
\delta(z,w)= \frac{(1-x^2-y^2)^2}{4}\, \delta(x-u)\,\delta(y-v)
\end{equation}
where $\delta(x)$ is the usual Dirac distribution 
with respect to Lebesgue measure on $\RR$.
In polar coordinates, where $w$ is taken as the origin,
$\tau=d(z,w)$, and
\begin{equation}
\delta(z,w)=\frac{\delta(\tau)}{2\pi \sinh\tau} .
\end{equation}
Property (G2) therefore says that $(\Delta+\lambda) G(z,w;\lambda)=0$
for $z\neq w$, and
\begin{equation}\label{greendiff2}
\int_{d(z,w)<\epsilon}
(\Delta+\lambda) G(z,w;\lambda) d\mu(z)= 1 \qquad\forall\epsilon>0.
\end{equation}
In view of (G3), there is a function $f\in\C^\infty(\RR_{>0})$ 
such that $f(\tau)=G(z,w;\lambda)$. Then
\begin{align}
1& =\int_{d(z,w)<\epsilon} (\Delta+\lambda) G(z,w;\lambda) d\mu(z) \\
& = 
2\pi \int_0^{\epsilon}
\left(\frac{d}{d\tau}
\left(\sinh\tau  \frac{d}{d\tau}\right) +\lambda\sinh\tau\right) 
f(\tau) \,d\tau  \\
& = 2\pi \sinh\epsilon f'(\epsilon) 
+ 2\pi \lambda \int_0^\epsilon \sinh\tau\,f(\tau)\,d\tau.
\end{align}
Taylor expansion around $\epsilon=0$ yields
$f'(\epsilon)=1/(2\pi\epsilon) + O(1)$
and thus $f(\epsilon)= (1/2\pi)\log \epsilon + O(1)$ as $\epsilon\to 0$.
Equation (G2) is therefore equivalent to
\begin{equation} \label{green}
\begin{cases}
(\Delta+\lambda) G(z,w;\lambda) = 0 , & d(z,w)>0 ,\\
G(z,w;\lambda) = (1/2\pi) \log d(z,w) +O(1), & d(z,w)\to 0 .
\end{cases}
\end{equation}

\begin{proposition}
If $\rho\in\CC$ with $\Im \rho < 1/2$, and $\lambda=\rho^2+\tfrac14$, then 
\begin{equation}
G(z,w;\lambda) = -\frac{1}{2\pi} Q_{-\frac12+\i\rho}\big(\cosh d(z,w)\big)
\end{equation} 
satisfies {\rm (G1)-(G4)}, where $Q_\nu$ is the Legendre function
of the second kind. 
\end{proposition}

\begin{proof}
With $f(\tau)=G(z,w;\lambda)$, (\ref{green}) becomes 
\begin{equation}
\left[\frac{1}{\sinh\tau}\frac{d}{d\tau}
\left(\sinh\tau  \frac{d}{d\tau}\right) 
+\lambda\right] f(\tau) =0 .
\end{equation}
Set $r=\cosh\tau$, $\tilde f(\cosh\tau)= f(\tau)$, and 
$\lambda=-\nu(\nu+1)$, to 
obtain Legendre's differential equation
\begin{equation}
\left[(1-r^2)\frac{d^2}{dr^2}
-2r \frac{d}{dr} + \nu(\nu+1) \right]
\tilde f(r) =  0 ,
\end{equation}
whose solutions are the associated Legendre functions $P_\nu(r)$
and $Q_\nu(r)$. $Q_\nu$ has the integral representation
\begin{equation}\label{intrep}
Q_{-\frac12+\i\rho}(\cosh\tau)
= \frac{1}{\sqrt 2} \int_\tau^\infty 
\frac{\e^{-\i \rho t}}{\sqrt{\cosh t-\cosh\tau}}\, dt
\end{equation}
which converges absolutely for $\Im \rho < 1/2$ and $\tau>0$. 
From this it is evident that $Q_{-\frac12+\i\rho}(\cosh\tau)\to 0$
as $\tau\to\infty$ (see also Lemma \ref{ubound} below), so (G4) holds.
For $t\to 0$ ($\rho$ fixed) it has the asymptotics required in \eqref{green},
cf.~the well known relation
\begin{equation}\label{asymp0}
Q_{-\frac12+\i\rho}(\cosh\tau) \sim 
- \big(\log(\tau/2) + \gamma + \psi(\tfrac12+\i\rho) \big)
\end{equation}
where $\gamma$ is Euler's constant and $\psi$ the logarithmic derivative
of Euler's $\Gamma$ function. 
\qed
\end{proof}

\begin{lemma}\label{ubound}
Given any $\epsilon>0$ there is a constant $C_\epsilon>0$
such that 
\begin{equation}
|Q_{-\frac12+\i\rho}(\cosh\tau)|
\leq C_\epsilon 
(1+|\log\tau|) \e^{\tau(\Im\rho-\tfrac12+\epsilon)} 
\end{equation}
uniformly for all $\tau> 0$ and $\rho\in\CC$ with $\Im\rho<\frac12-\epsilon$.
\end{lemma}

\begin{proof}
From the integral representation  \eqref{intrep} we infer
\begin{align}\label{intrep2}
|Q_{-\frac12+\i\rho}(\cosh\tau)|
& \leq \frac{1}{\sqrt 2} \int_\tau^\infty 
\frac{\e^{t\Im\rho}}{\sqrt{\cosh t-\cosh\tau}}\, dt \\
& \leq \frac{1}{\sqrt 2}\, 
\e^{\tau(\Im\rho-\tfrac12+\epsilon)} \int_\tau^\infty 
\frac{\e^{t(\tfrac12-\epsilon)}}{\sqrt{\cosh t-\cosh\tau}}\, dt 
\end{align}
since
\begin{equation}
\e^{t(\Im\rho-\tfrac12+\epsilon)}
\leq \e^{\tau(\Im\rho-\tfrac12+\epsilon)}
\end{equation}
for all $t\geq\tau$, if $\Im\rho-\tfrac12+\epsilon<0$ as assumed. 
The remaining integral
\begin{equation}
\int_\tau^\infty 
\frac{\e^{t(\tfrac12-\epsilon)}}{\sqrt{\cosh t-\cosh\tau}}\, dt 
\end{equation}
has a $\log\tau$ singularity at $\tau=0$ and is otherwise uniformly bounded
for all $\tau\to\infty$.
\qed
\end{proof}

\begin{remark}
This is only a crude upper bound, but sufficient for our purposes.
\end{remark}

To highlight the $\rho$ dependence, we shall use 
\begin{equation}
G_\rho(z,w)
= -\frac{1}{2\pi} Q_{-\frac12+\i\rho}\big(\cosh d(z,w)\big) 
\end{equation}
instead of $G(z,w;\lambda)$.

\begin{lemma}\label{Gfbound}
Suppose $f:\HH^2\to\CC$ with $|f(z)| \leq A \e^{\alpha d(z,o)}$,
with constants $A,\alpha>0$. Then the integral
\begin{equation}\label{this}
\int_{\HH^2} G_\rho(z,w) f(w)\, d\mu(w) 
\end{equation}
converges absolutely, 
uniformly in $\Re\rho$, provided $\Im\rho<-(\alpha+\tfrac12)$.
The convergence is also uniform in $z$ in compact sets in $\HH^2$.
\end{lemma}

\begin{proof}
Note that $|f(w)| \leq A \e^{\alpha d(o,w)}
\leq A \e^{\alpha d(o,z)} \e^{\alpha d(z,w)}$.
In polar coordinates (take $w$ as the origin) $\tau=d(z,w)$, and 
$d\mu=\sinh\tau\,d\tau\,d\phi$. In view of Lemma \ref{ubound},
the integral \eqref{this} is bounded by
\begin{equation}
\ll_\epsilon \int_0^\infty  |\log\tau|\;  
\e^{-(\tfrac12-\Im\rho-\epsilon-\alpha)\tau} \sinh\tau\,d\tau
\end{equation}
which converges under the hypothesis on $\Im\rho$.
\qed
\end{proof}

\section{Selberg's point-pair invariants \label{secSelberg}}

Let $H$ be a subgroup of $\Isom(\HH^2)$. An $H$-point-pair invariant 
$k: \HH^2\times\HH^2 \to \CC$ is defined by the relations
\begin{description}
\item[(K1)]
$k(g z,g w)=k(z,w)$ for all $g\in H$, $z,w\in\HH^2$;
\item[(K2)]
$k(w,z)=k(z,w)$ for all $z,w\in\HH^2$.
\end{description}
Here we will only consider point-pair invariants which are
functions of the distance between $z,w$, such as 
the Green's function $G_\rho(z,w)$ studied in the previous section.
Hence $H=\Isom(\HH^2)$ in this case. We sometimes use the notation
$k(\tau)=k(z,w)$ with $\tau=d(z,w)$.
Let us consider
\begin{equation}\label{defpp}
k(z,w) = \frac{1}{\pi\i} \int_{-\infty}^{\infty} 
G_\rho(z,w) \rho\, h(\rho)\, d\rho
\end{equation}
where the test function $h:\CC\to\CC$ satisfies the following conditions.
\begin{description}
\item[(H1)]
$h$ is analytic for $|\Im\rho|\leq \sigma$ for some $\sigma>1/2$;
\item[(H2)] 
$h$ is even, i.e., $h(-\rho)=h(\rho)$;
\item[(H3)]
$|h(\rho)|\ll (1+|\Re\rho|)^{-2-\delta}$ for some fixed $\delta>0$,
uniformly for all $\rho$ in the strip $|\Im\rho|\leq \sigma$.
\end{description}
For technical reasons we will sometimes use the stronger hypothesis
\begin{description}
\item[(H3*)]
$|h(\rho)|\ll_N (1+|\Re\rho|)^{-N}$ for any fixed $N>1$,
uniformly for all $\rho$ in the strip $|\Im\rho|\leq \sigma$.
\end{description}

The Fourier transform of $h$ is
\begin{equation}
g(t)=\frac{1}{2\pi} \int_\RR h(\rho)\, \e^{-\i\rho t} d\rho .
\end{equation}
With the integral representation \eqref{intrep} one immediately finds
\begin{equation}\label{kg}
k(z,w) =-\frac{1}{\pi\sqrt 2} \int_\tau^\infty 
\frac{g'(t)}{\sqrt{\cosh t-\cosh\tau}}\, dt , \quad \tau=d(z,w).
\end{equation}
The analyticity of $h$ and (H3) imply that $g$ and its first derivative 
({\em all derivatives} provided (H3*) holds)
are exponentially decaying for $|t|\to\infty$.
To see this, consider
\begin{align}
g^{(\nu)}(t)
&=\frac{1}{2\pi} \int_\RR (-\i\rho)^\nu h(\rho)\, \e^{-\i\rho t} d\rho \\
&=\frac{1}{2\pi} \int_{\RR-\i\sigma} 
(-\i\rho)^\nu h(\rho)\, \e^{-\i\rho t} d\rho \\
&= \frac{1}{2\pi} \e^{-\sigma t} 
\int_\RR [-\i(\rho-\i \sigma)]^\nu  h(\rho-\i \sigma)
\, \e^{-\i\rho t} d\rho .
\end{align}
Since, due to (H3*),
\begin{equation}
\int_\RR \left|(\rho-\i \sigma)^\nu  h(\rho-\i \sigma)\right|
d\rho <\infty
\end{equation}
we find
\begin{equation}\label{gbound}
|g^{(\nu)}(t)| \ll_\nu \e^{-\sigma |t|} .
\end{equation}

The point-pair invariant $k(z,w)$ gives rise to the linear operator
$L$ defined by
\begin{equation}
[Lf](z):=\int_{\HH^2} k(z,w) f(w) d\mu(w) .
\end{equation}

\begin{proposition}\label{propL1}
Suppose $f\in\C^2(\HH^2)$ is a solution of $(\Delta+\rho^2+\tfrac14) f=0$
with $|\Im\rho|\leq \sigma$ and $|f(z)|\leq A \e^{\alpha d(z,o)}$,
with constants $A>0,0\leq\alpha< \sigma-\tfrac12$.
Then, for $h$ satisfying {\rm (H1), (H2), (H3)},
\begin{equation}
Lf = h(\rho) f . 
\end{equation}
\end{proposition}

\begin{proof}
We have
\begin{align}
[Lf](z)& =\int_{\HH^2} k(z,w) f(w) d\mu(w) \\
& = \frac{1}{\pi\i} \int_{\HH^2} 
\bigg\{ \int_{-\infty}^{\infty} 
G_{\rho'}(z,w) \rho'\, h(\rho')\, d\rho' \bigg\} f(w) d\mu(w)\\
& = \frac{1}{\pi\i} \int_{\HH^2} 
\bigg\{ \int_{-\infty-\i \sigma}^{\infty-\i \sigma} 
G_{\rho'}(z,w) \rho'\, h(\rho')\, d\rho' \bigg\} f(w) d\mu(w)
\end{align}
where we 
have shifted the contour of integration by $\sigma$. Then
\begin{equation}
= \frac{1}{\pi\i} \int_{-\infty-\i \sigma}^{\infty-\i \sigma}  
\bigg\{ \int_{\HH^2} G_{\rho'}(z,w) f(w) d\mu(w) \bigg\} \rho'\, h(\rho')\, d\rho'
\end{equation}
since the inner integral converges absolutely, uniformly in $\Re\rho$,
see Lemma \ref{Gfbound}.
We have 
\begin{equation}
\int_{\HH^2} G_{\rho'}(z,w) f(w) d\mu(w) = (\Delta+{\rho'}^2+\tfrac14)^{-1}f(z)
= ({\rho'}^2-\rho^2)^{-1} f(z) 
\end{equation}
and thus
\begin{equation}
[Lf](z) = \frac{1}{\pi\i}\, f(z)\, \int_{-\infty-\i \sigma}^{\infty-\i \sigma} 
\frac{\rho' h(\rho')}{{\rho'}^2-\rho^2} \,d\rho'.
\end{equation}
This integral converges absolutely, cf.~(H3), and is easily calculated. 
We shift the contour from $-\i\sigma$ to $+\i\sigma$ 
and collect residues, so that
\begin{equation}
\frac{1}{2\pi\i} \int_{-\infty-\i \sigma}^{\infty-\i \sigma} 
\frac{2 \rho' h(\rho')}{{\rho'}^2-\rho^2} \,d\rho'
= h(\rho)+h(-\rho)
+\frac{1}{2\pi\i} \int_{-\infty+\i \sigma}^{\infty+\i \sigma} 
\frac{2 \rho' h(\rho')}{{\rho'}^2-\rho^2}  \,d\rho'.
\end{equation}
Since $h$ is even, the integral on the right hand side equals the 
negative of the left hand side, and thus
\begin{equation}
\frac{1}{2\pi\i} \int_{-\infty-\i \sigma}^{\infty-\i \sigma} 
\frac{2 \rho' h(\rho')}{{\rho'}^2-\rho^2}
= h(\rho) ,
\end{equation}
which concludes the proof. 
\qed
\end{proof}

It is useful to define the auxiliary functions $\Phi,Q:\RR_{\geq 0}\to\CC$
by the relations
\begin{equation}
\Phi\big(2\big(\cosh\tau-1)\big) =k(\tau)=k(z,w), \qquad \tau=d(z,w),
\end{equation}
and
\begin{equation}\label{Qdef}
Q\big(2(\cosh t-1)\big) = g(t) .
\end{equation}

\begin{lemma} \label{hQ}
The following statements are equivalent.
\begin{enumerate}
\item[{\rm (i)}]
$h$ satisfies {\rm (H1), (H2), (H3*)}.
\item[{\rm (ii)}]
$Q\in\C^\infty(\RR_{\geq 0})$ with
\begin{equation}\label{upperq}
\big|Q^{(\nu)}(\eta)\big| \ll_\nu \big(1+\eta\big)^{-\sigma-\nu} 
\qquad \forall \eta\geq 0 .
\end{equation}
\end{enumerate}
\end{lemma}

\begin{proof}
Clearly $g\in\C^\infty(\RR)$
if and only if $Q\in\C^\infty(\RR_{\geq 0})$
(this is obvious for $t\neq 0$; the problem at $t=0$ can be 
resolved by expanding in Taylor series).
In view of \eqref{gbound}, the bound \eqref{upperq} is evident for $\nu=0$.
The $\nu$th derivative of $g$ is of the form
\begin{equation}\label{gdiff}
g^{(\nu)}(t) = \sum_{j=0}^\nu a_{j\nu} \, \e^{j |t|} \, Q^{(j)}(2(\cosh t-1)) 
\;(1+O(\e^{-|t|}))
\end{equation}
with suitable coefficients $a_{j\nu}$. Hence, by induction on $\nu$,
\begin{align}
\e^{\nu |t|} \, \big|Q^{(\nu)}(2(\cosh t-1)) \big|
& \ll_\nu \big| g^{(\nu)}(t) \big| + 
\big|\sum_{j=0}^{\nu-1} a_{j\nu} 
\, \e^{j |t|} \, Q^{(j)}(2(\cosh t-1))\big| \\
& \ll_\nu \big| g^{(\nu)}(t) \big| + \sum_{j=0}^{\nu-1} \e^{j |t|} 
\e^{(-\sigma-j)|t|}  \\
& \ll_\nu \e^{-\sigma |t|} .
\end{align}
This proves (i) $\Rightarrow$ (ii). Conversely, \eqref{upperq}
implies via \eqref{gdiff} the exponential decay of $g$, which proves
(H1). (H3*) follows from $g\in\C^\infty(\RR)$.
\qed
\end{proof}

The integral transform \eqref{kg} reads in terms of 
the functions $Q,\Phi$,
\begin{equation}\label{inttrans}
\Phi(\xi) = -\frac{1}{\pi} \int_\xi^\infty 
\frac{Q'(\eta)}{\sqrt{\eta-\xi}}\, d\eta .
\end{equation}
In order to specify $Q$ uniquely for a given $\Phi$, we always assume in the following that $Q(\eta) \to 0$ for $\eta\to\infty$.

\begin{lemma} \label{QPsi}
Consider the following conditions.
\begin{enumerate}
\item[{\rm (i)}]
$Q\in\C^\infty(\RR_{\geq 0})$ and 
$\big|Q^{(\nu)}(\eta)\big| \ll_\nu \big(1+\eta\big)^{-\sigma-\nu}$.
\item[{\rm (ii)}] 
$\Phi\in\C^\infty(\RR_{\geq 0})$ and 
$\big|\Phi^{(\nu)}(\xi)\big| 
\ll_\nu \big(1+\xi\big)^{-\sigma-\nu-1/2+\epsilon}$.
\end{enumerate}
Then {\rm (i)} implies {\rm (ii)} for any fixed $\epsilon>0$,
and {\rm (ii)} implies {\rm (i)} for any fixed $\epsilon<0$.
\end{lemma}

\begin{proof}
The $\nu$th derivative of $\Phi$ is
\begin{align}
\Phi^{(\nu)}(\xi) & = -\frac{1}{\pi} \frac{d^\nu}{d\xi^\nu} \int_0^\infty 
\frac{Q'(\eta+\xi)}{\sqrt{\eta}}\, d\eta \\
& = -\frac{1}{\pi} \int_0^\infty 
\frac{Q^{(\nu+1)}(\eta+\xi)}{\sqrt{\eta}}\, d\eta \label{Q2} \\
& = -\frac{1}{\pi} \int_\xi^\infty 
\frac{Q^{(\nu+1)}(\eta)}{\sqrt{\eta-\xi}}\, d\eta .
\end{align}
Therefore (i) implies $\Phi\in\C^\infty(\RR_{\geq 0})$.
Furthermore, from \eqref{Q2},
\begin{align}
\big|\Phi^{(\nu)}(\xi)\big|
& \ll_\nu \int_0^\infty 
\frac{\big(1+\eta+\xi\big)^{-\sigma-\nu-1}}{\sqrt{\eta}}\, d\eta \\
& \ll_\nu 
\big(1+\xi\big)^{-\sigma-\nu-1/2+\epsilon}
\int_0^\infty 
\frac{\big(1+\eta\big)^{-1/2-\epsilon}}{\sqrt{\eta}}\, d\eta ,
\end{align}
for $\epsilon>0$ small enough. The last integral converges for any 
$\epsilon>0$.

The implication {\rm (ii)} $\Rightarrow$ {\rm (i)}
follows analogously from the inversion formula 
\begin{equation}\label{inv}
Q(\eta) = \int_\eta^\infty 
\frac{\Phi(\xi)}{\sqrt{\xi-\eta}}\, d\xi .
\end{equation}
To show that \eqref{inv} is indeed consistent with \eqref{inttrans},
write \eqref{inv} in the form
\begin{equation}\label{QPhi}
Q(\eta)=\int_{-\infty}^\infty \Phi(\eta+\xi^2) d\xi
\end{equation}
and thus
\begin{equation}\label{Qprime}
Q'(\eta)=\int_{-\infty}^\infty \Phi'(\eta+\xi^2) d\xi .
\end{equation}
The right hand side of \eqref{inttrans} is 
\begin{align}
-\frac{1}{\pi} \int_\xi^\infty 
\frac{Q'(\eta)}{\sqrt{\eta-\xi}}\, d\eta 
&=-\frac{1}{\pi}
\int_\RR Q'(\xi+\eta^2) d\eta \\
&=-\frac{1}{\pi} 
\int_\RR \int_\RR \Phi'(\xi+\eta_1^2+\eta_2^2) d\eta_1 d\eta_2 ,
\end{align}
where we have used \eqref{Qprime} in the last step.
This equals of course
\begin{equation}
= - 2 \int_0^\infty \Phi'(\xi+r^2) r dr 
=- \int_0^\infty \Phi'(\xi+r) dr = \Phi(\xi)
\end{equation}
which yields the left hand side of \eqref{inttrans}.
\qed
\end{proof}

\begin{proposition} \label{hk}
The following statements are equivalent.
\begin{enumerate}
\item[{\rm (i)}]
$h$ satisfies {\rm (H1), (H2), (H3*)}.
\item[{\rm (ii)}]
$k(z,w)$ is in $\C^\infty(\HH^2\times\HH^2)$ with the bound on the $\nu$th 
derivative,
\begin{equation}\label{upperk}
\big|k^{(\nu)}(\tau)\big| \ll_\nu \e^{-(\sigma+1/2-\epsilon)\tau} 
\qquad \forall \tau\geq 0 ,
\end{equation}
for any fixed $\epsilon>0$.
\end{enumerate}
\end{proposition}

\begin{proof}
In view of Lemmas \ref{hQ} and \ref{QPsi}, the statement (i) is equivalent
to the condition for $\Phi$, statement (ii) in Lemma \ref{QPsi}.
Since $k(z,w)=k(\tau)=\Phi\big(2(\cosh\tau -1)\big)$, the proof is exactly the
same as that of Lemma \ref{hQ} with $g(t)$ replaced by $k(\tau)$,
and $Q\big(2(\cosh t -1)\big)$ by $\Phi\big(2(\cosh\tau-1)\big)$.
\qed
\end{proof}

\section{The ghost of the sphere \label{ghost}}

Note that for $z=w$, the kernel
$k(z,w)$ has a finite value, unlike the logarithmic divergence
of the Green's function $G_\rho(z,w)$. In fact,
\begin{align}
k(z,z) & =-\frac{1}{\pi\sqrt 2} \int_0^\infty 
\frac{g'(t)}{\sqrt{\cosh t-1}}\, dt \\
& =-\frac{1}{2\pi} \int_0^\infty 
\frac{g'(t)}{\sinh(t/2)}\, dt \\ \label{ditte}
& = \frac{1}{4\pi^2} \int_0^\infty \bigg\{ \int_{-\infty}^\infty 
\frac{\sin(\rho t)}{\sinh(t/2)} h(\rho)\, \rho\,d\rho\bigg\} dt \\
& = \frac{1}{4\pi^2} \int_{-\infty}^\infty  \bigg\{ \int_0^\infty
\frac{\sin(\rho t)}{\sinh(t/2)}\, dt \bigg\}  h(\rho)\, \rho\,d\rho
\end{align}
where changing the order of integration is justified, since 
\begin{equation}
\int_0^\infty
\left| \frac{\sin(\rho t)}{\sinh(t/2)} \right| dt
\leq  |\rho| \int_0^\infty
\frac{t}{\sinh(t/2)}\, dt \ll |\rho|
\end{equation}
and $|h(\rho)|\ll(1+|\rho|)^{-4}$, assuming (H3*). 
We use the geometric series expansion
\begin{equation}\label{sinh}
\frac{1}{\sinh(t/2)}
=
\frac{2\e^{-t/2}}{1-\e^{-t}}
=
2 \sum_{l=0}^\infty \exp\left[-\left(l+\tfrac12\right)t\right],
\end{equation}
so for $|\Im\rho|<1/2$
\begin{align}
\int_0^\infty
\frac{\sin(\rho t)}{\sinh(t/2)}\, dt 
& = 
\i \sum_{l=0}^\infty  
\left[\frac{1}{\i\rho-(l+\tfrac12)} + \frac{1}{\i\rho+(l+\tfrac12)}\right]
\label{here} \\
& = -\pi\i \tan(\pi\i\rho) = \pi \tanh(\pi\rho) ,
\end{align}
compare \eqref{here14}.
We conclude
\begin{equation}\label{weyl}
k(z,z)=
\frac{1}{4\pi} \int_{-\infty}^\infty h(\rho) \tanh(\pi\rho)\,\rho\, d\rho .
\end{equation}

Let us conclude this section by noting that the logarithmic 
divergence of the Green's function is independent of $\rho$,
see \eqref{green}.
It may therefore be removed by using instead
\begin{equation}
G_\rho(z,w)- G_{\rho_*}(z,w)
\end{equation}
where $\rho_*\neq \rho$ is a fixed constant in $\CC$ with $|\Im\rho_*|<1/2$.
We then have from \eqref{intrep}
\begin{align} 
\lim_{w\to z} \left[ G_\rho(z,w)- G_{\rho_*}(z,w) \right]
& =
-\frac{1}{2\pi \sqrt 2} \int_\tau^\infty 
\frac{\e^{-\i \rho t}-\e^{-\i \rho_* t}}{\sqrt{\cosh t-\cosh\tau}}\, dt \\
&=-\frac{1}{4\pi}\int_0^\infty
\frac{\e^{-\i \rho t}-\e^{-\i \rho_* t}}{\sinh(t/2)}\, dt \\
& = 
- \frac{1}{2\pi\i} \sum_{l=0}^\infty  
\left[\frac{1}{\rho-\i(l+\tfrac12)} - \frac{1}{\rho_*-\i(l+\tfrac12)}\right],
\label{ztow}
\end{align}
where we have used the geometric series expansion \eqref{sinh}
as above. The last sum clearly converges since
\begin{equation}
\frac{1}{\rho-\i(l+\tfrac12)} - \frac{1}{\rho_*-\i(l+\tfrac12)} =O(l^{-2}).
\end{equation}

\begin{remark}
In analogy with the trace formula for the sphere, we may
view the geometric series expansion, 
\begin{equation}\label{tanh}
\tanh(\pi\rho)  
= \sum_{n\in\ZZ} (-1)^n \e^{-2\pi |n| \rho} ,
\end{equation}
cf.~\eqref{tan}, as a sum
over closed orbits on the sphere, but now with imaginary action.
These orbits have an interpretation as tunneling (or ghost) orbits.
\end{remark}

\section{Hyperbolic surfaces \label{secHyperbolic2}}

Let $\scrM$ be a smooth Riemannian surface (finite or infinite)  
of constant negative curvature which can be represented
as the quotient $\GamHH$, where $\Gamma$ is a strictly hyperbolic 
Fuchsian group (i.e., all elements $\gamma\in\Gamma-\{1\}$ have
$\ell_\gamma>0$).
The space of square integrable functions on $\scrM$
may therefore be identified with the space of measurable functions
$f:\HH^2\to\CC$ satisfying the properties
\begin{equation}\label{comm}
T_\gamma f = f \quad\forall\gamma\in\Gamma
\end{equation}
and 
\begin{equation}
\|f\|^2:=\int_{\FGam} |f|^2 d\mu < \infty
\end{equation}
where $T_\gamma$ is the translation operator defined in \eqref{trans} and
$\FGam$ is any fundamental domain of $\Gamma$ in $\HH^2$. We denote
this space by $\L^2(\GamHH)$. The inner product
\begin{equation}
\langle f_1,f_2 \rangle = \int_{\FGam} f_1 \overline f_2 d\mu 
\end{equation}
makes $\L^2(\GamHH)$ a Hilbert space.
Similarly, we may identify $\C^\infty(\GamHH)$ with the space of functions
$f\in\C^\infty(\HH^2)$ satisfying \eqref{comm} (note that more care has to be
taken here when $\Gamma$ contains elliptic elements).
Since $\Delta$ commutes with $T_\gamma$, it maps 
$\C^\infty(\GamHH)\to\C^{\infty}(\GamHH)$. 

To study the spectrum of the Laplacian on $\GamHH$, let us consider
the linear operator $L$ of functions on $\GamHH$,
\begin{equation}
[L f](z):=\int_{\GamHH} k_\Gamma(z,w) f(w) d\mu(w) 
\end{equation}
with kernel
\begin{equation}\label{defk}
k_\Gamma(z,w) = \sum_{\gamma\in\Gamma} k(\gamma z, w) 
\end{equation}
with the point-pair invariant $k$ as defined in \eqref{defpp}.
The convergence of the sum is guaranteed by the following lemma, 
cf.~Proposition \ref{convert} below.

\begin{lemma}\label{uppi}
For every $\delta>0$, there is a $C_{\delta}>0$ 
such that
\begin{equation}
\sum_{\gamma\in\Gamma} \e^{-(1+\delta) d(\gamma z, w)} 
\leq C_{\delta}
\end{equation}
for all $(z,w)\in\HH^2\times\HH^2$.
\end{lemma}

\begin{proof}
Place a disk $\scrD_\gamma(r)=\{ z'\in\HH: d(\gamma z,z')\leq r \}$ around 
every point $z_\gamma=\gamma z$, and denote the area of $\scrD_\gamma(r)$
by $\Area(r)$. 
Then
\begin{equation}
\e^{-(1+\delta) d(\gamma z, w)} 
\leq 
\frac{\e^{r(1+\delta)}}{\Area(r)}
\int_{\scrD_\gamma(r)} \e^{-(1+\delta) d(z', w)} d\mu(z') .
\end{equation}
If 
\begin{equation}
r<\frac12 \min_{\gamma\in\Gamma-\{1\}} \ell_\gamma
\end{equation} 
the disks $\scrD_\gamma(r)$ do not overlap.
(Note that $\min_{\gamma\in\Gamma-\{1\}} \ell_\gamma>0$
since $\Gamma$ is strictly hyperbolic
and acts properly discontinuously on $\HH^2$.)
Therefore
\begin{align}
\sum_{\gamma\in\Gamma} \e^{-(1+\delta) d(\gamma z, w)} & \leq 
\frac{\e^{r(1+\delta)}}{\Area(r)} 
\int_{\HH^2} \e^{-(1+\delta) d(z', w)} d\mu(z') \\
& = \frac{2\pi\e^{r(1+\delta)}}{\Area(r)}  \int_0^\infty \e^{-(1+\delta) \tau}
\sinh\tau \, d\tau .
\end{align}
This integral converges for any $\delta>0$. \qed
\end{proof}

\begin{proposition}\label{convert}
If $h$ satisfies {\rm (H1), (H2), (H3*)}, then
the kernel $k_\Gamma(z,w)$ is in $\C^\infty(\GamHH \times\GamHH)$,
with $k_\Gamma(z,w)=k_\Gamma(w,z)$.
\end{proposition}

\begin{proof}
Proposition \ref{hk} and Lemma \ref{uppi} show that the sum over
$k(\gamma z,w)$ converges absolutely and uniformly
(take $\delta=\sigma-1/2-\epsilon>0$).
The same holds for sums over any derivative of $k(\gamma z,w)$.
Hence $k_\Gamma(z,w)$ is in $\C^\infty(\HH^2\times\HH^2)$.
To prove invariance under $\Gamma$, note that
\begin{equation}
k_\Gamma(\gamma z,w) = \sum_{\gamma'\in\Gamma} k(\gamma'\gamma z, w) 
= \sum_{\gamma'\in\Gamma} k(\gamma' z, w) =  k_\Gamma(z,w) .
\end{equation}
Thus $k_\Gamma(z,w)$ is a function on $\GamHH$ with respect to the 
first argument.
Secondly
\begin{multline}
k_\Gamma(z,w) = \sum_{\gamma'\in\Gamma} k(\gamma'z, w) 
= \sum_{\gamma'\in\Gamma} k(w,\gamma' z) \\ 
= \sum_{\gamma'\in\Gamma} k({\gamma'}^{-1}w,{\gamma'}^{-1}\gamma' z) 
= \sum_{\gamma'\in\Gamma} k({\gamma'}^{-1} w, z) 
=k_\Gamma(w,z) ,
\end{multline}
which proves symmetry. Both relations imply immediately
$k_\Gamma(z,\gamma w) = k_\Gamma(z,w)$.
\qed
\end{proof}

\begin{proposition}\label{propL11}
Suppose $f\in\C^2(\GamHH)$ is a solution of $(\Delta+\rho^2+\tfrac14) f=0$
with $|\Im\rho|\leq \sigma$ and $|f(z)|\leq A \e^{\alpha d(z,o)}$,
with constants $A>0,0\leq\alpha< \sigma-\tfrac12$.
Then, for $h$ satisfying {\rm (H1), (H2), (H3)},
\begin{equation}
Lf = h(\rho) f . 
\end{equation}
\end{proposition}

\begin{proof}
Note that 
\begin{equation}
\int_{\FGam} k_\Gamma(z,w) f(w) d\mu(w) =
\int_{\HH^2} k(z,w) f(w) d\mu(w) 
\end{equation}
and recall Proposition \ref{propL1}.
\qed
\end{proof}

\section{A trace formula for hyperbolic cylinders \label{secCylinders}}

The simplest non-trivial example of a hyperbolic surface is
a hyperbolic cylinder. To construct one,
fix some $\gamma\in\Isom^+(\HH^2)$ of length $\ell=\ell(\gamma)>0$, 
and set $\Gamma=\Z$,
where $\Z$ is the discrete subgroup generated by $\gamma$, i.e.,
\begin{equation}
\Z= \{ \gamma^n :\; n\in\ZZ \} .
\end{equation}
We may represent $\ZHH$ in halfplane coordinates, which are chosen
in such a way that 
\begin{equation}
\gamma=
\begin{pmatrix}
\e^{\ell/2} & 0 \\ 0 & \e^{-\ell/2} 
\end{pmatrix}.
\end{equation}
A fundamental domain for the action
of $\gamma$ on $\H$, $z\mapsto \e^\ell z$, is given by
\begin{equation}
\{ z\in\H: \; 1 \leq y < \e^\ell \} .
\end{equation}
It is therefore evident that $\ZHH$ has infinite volume. 
A more convenient set of parameters for the cylinder 
are the coordinates $(s,u)\in\RR^2$, with
\begin{equation}\label{su}
x=u \e^s, \qquad y=\e^s ,
\end{equation}
where the volume element reads now
\begin{equation}
d\mu= ds\,du.
\end{equation}
In these coordinates, the action of $\gamma$ is $(s,u)\mapsto(s+\ell,u)$,
and hence a fundamental domain is
\begin{equation}
\scrF_{\Z}=\{ (s,u)\in\RR^2 : \; 0\leq s < \ell \} .
\end{equation}
Note that
\begin{equation}
\cosh d(\gamma^n z,z) = 1 + \frac{|\e^{n\ell} z-z|^2}{2\e^{n\ell} y^2}
=
1+ 2\sinh^2(n\ell/2) (1+u^2) 
\end{equation}
and hence
\begin{align}\label{cylk}
k_\Z(z,z)& = k(z,z) + \sum_{n\neq 0}^\infty k(\gamma^n z, z) \\
& = k(z,z) + 2 \sum_{n=1}^\infty k(\gamma^n z, z) \\
& = k(z,z)+2 \sum_{n=1}^\infty \Phi\big(4\sinh^2(n\ell/2)(1+u^2) \big).
\end{align}
From this we can easily work out a trace formula for the
hyperbolic cylinder:
\begin{proposition}\label{propfirst}
If $h$ satisfies {\rm (H1), (H2), (H3)}, then
\begin{equation}\label{firsttrace}
\int_{\ZHH} \big[ k_\Z(z,z) - k(z,z) \big] d\mu
=\sum_{n=1}^\infty 
\frac{\ell\,g(n\ell)}{\sinh(n\ell/2)}.
\end{equation}
\end{proposition}

\begin{proof}
We have
\begin{align}
\int_\RR \int_0^\ell & \Phi\big(4\sinh^2(n\ell/2) (1+u^2) \big) \,ds\,du \\
& = 
\frac{\ell}{2\sinh(n\ell/2)}
 \int_\RR \Phi\big(4\sinh^2(n\ell/2) +\xi^2 \big) \,d\xi \\
\intertext{and with \eqref{QPhi},}
& = 
\frac{\ell}{2\sinh(n\ell/2)}\,
Q\big(4\sinh^2(n\ell/2) \big)\\
& = 
\frac{\ell}{2\sinh(n\ell/2)}\,
Q\big(2(\cosh(n\ell)-1) \big) 
\end{align}
which yields the right hand side of \eqref{firsttrace}, cf.~\eqref{Qdef}.
\qed
\end{proof}

\begin{proposition}\label{scatter}
If $h$ satisfies {\rm (H1), (H2), (H3)}, then
\begin{equation}
\int_{\ZHH} \big[ k_\Z(z,z) - k(z,z) \big] d\mu=
\int_{\RR}  h(\rho)\,n_\Z(\rho)\,d\rho
\end{equation}
where
\begin{equation}
n_\Z(\rho)=
\frac{\ell}{\pi}
\sum_{m=0}^\infty
\left\{\exp\left[\left( m+\tfrac12+\i\rho\right) 
\ell \right]-1\right\}^{-1}
\end{equation}
is a meromorphic function in $\CC$ with simple poles at the points
\begin{equation}
\rho_{\nu m}= \frac{\nu}{2\pi\ell}+\i\left(m+\tfrac12\right) ,
\qquad \nu\in\ZZ, \quad m=0,1,2,\ldots, 
\end{equation}
and residues $\res_{\rho_{\nu m}} n_\Z=1/(\pi\i)$.
\end{proposition}

\begin{proof}
The geometric series expansion of $1/\sinh$ \eqref{sinh} yields
\begin{equation} \label{nnn}
\sum_{n=1}^\infty 
\frac{\ell\,g(n\ell)}{\sinh(n\ell/2)}
=\frac{\ell}{\pi}
\int_\RR \sum_{m=0}^\infty \sum_{n=1}^\infty 
\exp\left[-\left(m+\tfrac12+\i\rho\right)n\ell \right]
h(\rho)\,d\rho
\end{equation}
and using again the geometric series, this time for the sum over $n$,
\begin{align}
\sum_{n=1}^\infty 
\exp\left[-\left(m+\tfrac12+\i\rho\right)n\ell \right]
& =
\left\{1-\exp\left[-\left(m+\tfrac12+\i\rho\right)\ell \right]
\right\}^{-1} - 1 \\
& =
\left\{\exp\left[\left(m+\tfrac12+\i\rho\right)\ell \right]-1
\right\}^{-1} .
\end{align}
This proves the formula for $n_\Z(\rho)$.
Near each pole $\rho_{\nu m}$ we have
\begin{align}
n_\Z(\rho) & \sim \frac{\ell}{\pi}
\left\{\exp\left[\left(m+\tfrac12+\i\rho\right)\ell \right] -1 
\right\}^{-1} \\
& \sim \frac{1}{\pi\i}\; \frac{1}{\rho-\left[(2\pi/\ell)\nu 
+\i\left(m+\tfrac12\right)\right]}
\end{align}
and so $\res_{\rho_{\nu m}} n_\Z=1/(\pi\i)$.
\qed
\end{proof}

The poles of $n_\Z(\rho)$ are called the {\em scattering poles}
of the hyperbolic cylinder.
A useful formula for $n_\Z$ is
\begin{equation} \label{nnn2}
n_\Z(\rho)=
\frac{1}{2\pi} \sum_{n=1}^\infty 
\frac{\ell\,\e^{-\i \rho n \ell}}{\sinh(n\ell/2)}, \qquad \Im\rho<1/2,
\end{equation}
which follows immediately from the above proof.
Furthermore, by shifting the path of 
integration to $-\i\infty$, we have the identity
\begin{equation} \label{nnn3}
\int_\RR \frac{n_\Z(\rho')}{\rho^2-{\rho'}^2}\, d\rho' =
\frac{\pi\i}{\rho} n_\Z(\rho),     \qquad \Im\rho<0.
\end{equation}

\section{Back to general hyperbolic surfaces \label{secBack}}

Let us now show that the kernel $k_\Gamma(z,w)$ of a general hyperbolic
surface $\GamHH$ (with $\Gamma$ strictly hyperbolic)
can be written as a superposition of kernels corresponding
to hyperbolic cylinders. 

Define the {\em conjugacy class} of any element $\gamma\in\Gamma$ as
\begin{equation}
\{\gamma\} := \{ \tilde\gamma \in\Gamma: \, \tilde\gamma=g\gamma g^{-1}
\text{ for some } g\in\Gamma \} .
\end{equation}
Clearly the length $\ell_\gamma$ is the same for all elements in one
conjugacy class.
The {\em centralizer} of $\gamma$ is
\begin{equation}
\Z_\gamma := \{ g \in\Gamma: \,  g \gamma=\gamma g \} .
\end{equation}

\begin{lemma}
If $\gamma\in\Gamma$ is hyperbolic, then the centralizer is the
infinite cyclic subgroup
\begin{equation}
\Z_\gamma = \{ \gamma_*^n :\; n\in\ZZ \} ,
\end{equation}
where $\gamma_*\in\Gamma$ is unique element such that
$\gamma_*^m=\gamma$ for some $m\in\NN$ and there is no $\tilde\gamma\in\Gamma$ such that
$\tilde\gamma^n=\gamma^*$ for any $n\in\NN$, $n>1$. 
\end{lemma}

(The element $\gamma_*$ is called the {\em primitive} of $\gamma$.)

\begin{proof}
$\gamma\in\PSL(2,\RR)\simeq\Isom^+(\HH^2)$ is conjugate to a diagonal matrix 
\begin{equation}
\begin{pmatrix}
\e^{\ell_\gamma/2} & 0 \\ 0 & \e^{-\ell_\gamma/2}
\end{pmatrix}
\end{equation}
with $\ell_\gamma>0$. The equation
\begin{equation}
\begin{pmatrix}
\e^{\ell_\gamma/2} & 0 \\ 0 & \e^{-\ell_\gamma/2}
\end{pmatrix}
\begin{pmatrix}
a & b \\ c & d
\end{pmatrix}
=
\begin{pmatrix}
a & b \\ c & d
\end{pmatrix}
\begin{pmatrix}
\e^{\ell_\gamma/2} & 0 \\ 0 & \e^{-\ell_\gamma/2}
\end{pmatrix}
\end{equation}
has the only solution $b=c=0$, $a=d^{-1}$. Hence the centralizer
is a diagonal subgroup of (a conjugate of) $\Gamma$. 
Since $\Gamma$ is discrete,
the centralizer must be discrete, which forces it to be cyclic.
\qed
\end{proof}

\begin{remark}
If $\gamma$ is hyperbolic and $\Gamma$ strictly hyperbolic, 
then $\{\gamma\}\neq \{\gamma^n\}$ for all $n\neq 1$. 
Furthermore, the centralizers of $\gamma$ and $\gamma^n$ coincide.
\end{remark}

The sum in \eqref{defk} can now be expressed as
\begin{align}
\sum_{\gamma\in\Gamma} k(\gamma z, w) 
& =
k(z, w) +
\sum_{\gamma\in H} \sum_{g\in \Z_\gamma\backslash\Gamma} 
k(g^{-1}\gamma g z, w)  \\
& =
k(z, w) +
\sum_{\gamma\in H} \sum_{g\in \Z_\gamma\backslash\Gamma} k(\gamma g z,g w) 
\label{kk}
\end{align}
where the respective first sums run over a set $H$ of hyperbolic elements,
which contains one representative for each conjugacy class $\{\gamma\}$.
We may replace this sum by a sum over {\em primitive} elements.
If we denote by $H_*\subset H$ the subset of primitive elements, \eqref{kk}
equals
\begin{align}
& = k(z, w)  +
\sum_{\gamma\in H} \sum_{g\in \Z_\gamma\backslash\Gamma} k(\gamma g z,g w) 
 \\ & =
k(z, w)  +
\sum_{\gamma\in H_*}
\sum_{g\in \Z_\gamma\backslash\Gamma}  \sum_{n=1}^\infty 
k(\gamma^n g z,g w) 
\end{align}
and hence, finally,
\begin{equation}\label{finale}
k_\Gamma(z, w) - k(z,w)=
\frac12\sum_{\gamma\in H_*}
\sum_{g\in \Z_\gamma\backslash\Gamma} 
\left\{k_{\Z_\gamma}(g z,g w) - k(g z,g w) \right\} ;
\end{equation}
recall that
\begin{equation}\label{finale2}
k_{\Z_\gamma}(z,w)=
\sum_{n\in\ZZ} k(\gamma^n z,w).
\end{equation}

\section{The spectrum of a compact surface \label{secSpectrum}}

It is well known that for any compact Riemannian manifold, 
$-\Delta$ has positive discrete spectrum, i.e.,
\begin{equation}
0=\lambda_0< \lambda_1 \leq \lambda_2 \leq \ldots \to \infty ,
\end{equation}
with corresponding eigenfunctions 
$\varphi_0=\text{const},
\varphi_1,\varphi_2,\ldots\in\C^\infty(\GamHH)$,
which satisfy
\begin{equation}
(\Delta+\lambda_j)\varphi_j=0 
\end{equation}
and form an orthonormal basis of $\L^2(\GamHH)$.
Furthermore, since $\Delta$ is
real-symmetric, the $\varphi_j$ can be chosen to be real-valued.
We furthermore define 
\begin{equation}
\rho_j=\sqrt{\lambda_j-\tfrac14}, \qquad -\pi/2 \leq \arg\rho_j < \pi/2.
\end{equation}
If $f\in\C^2(\GamHH)$, the expansion 
\begin{equation}\label{expand}
f(z)= \sum_j c_j \varphi_j(z), \qquad c_j=\langle f, \varphi_j \rangle,
\end{equation}
converges absolutely, uniformly for all $z\in\HH^2$.
This follows from general spectral theoretic arguments,
compare \cite[p.~3 and Chapter {\sc three}]{HejhalI}.

\begin{proposition} \label{propL2}
If $h$ satisfies {\rm (H1), (H2), (H3)}, then
\begin{equation}
L \varphi_j  = h(\rho_j) \varphi_j .
\end{equation}
\end{proposition}

\begin{proof}
Apply Proposition \ref{propL11}. Each eigenfunction $\varphi_j$ is
bounded so $\alpha=0$. Furthermore, by the positivity of $-\Delta$,
we have $|\Im\rho|\leq 1/2< \sigma$. 
\qed
\end{proof}

\begin{proposition}\label{pretrace1}
If $h$ satisfies {\rm (H1), (H2), (H3*)}, then 
\begin{equation}
k_\Gamma(z,w)
=\sum_{j=0}^\infty  h(\rho_j)\, \varphi_j(z)\,\overline\varphi_j(w),
\end{equation}
which converges absolutely, uniformly in $z,w\in\HH^2$.
\end{proposition}

\begin{proof}
The spectral expansion \eqref{expand} of $k_\Gamma(z,w)$ as a function
of $z$ yields
\begin{equation}
k_\Gamma(z,w)
=\sum_{j=0}^\infty c_j \varphi_j(z), 
\end{equation}
with
\begin{equation}
c_j=\int_{\GamHH} k_\Gamma(z,w)\overline\varphi_j(z) \,d\mu(z)
= \overline{[L \varphi_j]}(w) =\overline h(\rho_j)\,\overline\varphi_j(w)
= h(\rho_j)\,\overline\varphi_j(w).
\end{equation}
The proof of uniform convergence follows from standard spectral theoretic
arguments \cite[Prop. 3.4, p.12]{HejhalI}.
\qed
\end{proof}

In the case $z=w$, Proposition \ref{pretrace1} implies immediately
the following theorem.

\begin{theorem}[Selberg's pre-trace formula]\label{pretrace2}
If $h$ satisfies {\rm (H1), (H2), (H3*)}, then
\begin{equation}\label{pretraceEQ}
\sum_{j=0}^\infty  h(\rho_j)\, |\varphi_j(z)|^2
= \frac{1}{4\pi} \int_{-\infty}^\infty h(\rho) \tanh(\pi\rho)\,\rho\, d\rho
+ \sum_{\gamma\in\Gamma-\{1\}} k(\gamma z, z) .
\end{equation}
which converges absolutely, uniformly in $z\in\HH^2$.
\end{theorem}

\begin{proof}
Use \eqref{weyl} for the $\gamma=1$ term.
\qed
\end{proof}

Using \eqref{finale}, the pre-trace formula
\eqref{pretraceEQ} becomes
\begin{align}
\sum_{j=0}^\infty  h(\rho_j)\, |\varphi_j(z)|^2
& =\frac{1}{4\pi} \int_{-\infty}^\infty h(\rho) \tanh(\pi\rho)\,\rho\, d\rho\\
& \quad +\label{kZ}
\frac12 \sum_{\gamma\in H_*}
\sum_{g\in \Z_\gamma\backslash\Gamma}  
\left\{k_{\Z_\gamma}(g z,g z) - k(g z,g z) \right\}.
\end{align}

\begin{theorem}[Selberg's trace formula]\label{trace}
If $h$ satisfies {\rm(H1), (H2), (H3*)}, then
\begin{equation}\label{traceEQ}
\sum_{j=0}^\infty  h(\rho_j)
= \frac{\Area(\scrM)}{4\pi} 
\int_{-\infty}^\infty h(\rho) \tanh(\pi\rho)\,\rho\, d\rho
+ \sum_{\gamma\in H_*}  \sum_{n=1}^\infty
\frac{\ell_{\gamma}\,g(n\ell_{\gamma})}{2\sinh(n\ell_{\gamma}/2)} ,
\end{equation}
which converges absolutely.
\end{theorem}

(We will see in the next section (Corollary \ref{tracecor}) that
the condition (H3*) may in fact be replaced by (H3).)

\begin{proof}
We integrate both sides of the pre-trace formula \eqref{pretraceEQ}
over $\GamHH$. By the $\L^2$ normalization of the eigenfunctions $\varphi_j$, 
the left hand side of \eqref{pretraceEQ} yields the left hand side of
\eqref{traceEQ}. The first term on the right hand side is trivial,
and the second term follows from the observation that 
\begin{equation}
\sum_{g\in \Z_\gamma\backslash\Gamma} \int_{\GamHH}  f(g z)\, d\mu
= 
\int_{\ZHH} f(z)\, d\mu ,
\end{equation}
which allows us to apply Proposition \ref{propfirst} to 
the inner sum in \eqref{kZ}.
\qed
\end{proof}

\begin{remark}\label{rem}
The absolute convergence of the sum on the right hand side of
\eqref{traceEQ} only requires (H1), or
\begin{equation}
|g(t)| \ll \e^{-\sigma |t|} , \qquad \forall t > 0 .
\end{equation}
One way of seeing this is is that, 
since $g$ is only evaluated on the discrete subset
(the {\em length spectrum})
\begin{equation}
\{\ell_\gamma :\; \; \gamma\in\Gamma-\{1\} \}\subset\RR_{>0},
\end{equation} 
we may replace $g$ by an even $\C^\infty(\RR)$ function
$\tilde g$ (for which absolute convergence is granted)
so that $g(\ell_\gamma)=\tilde g(\ell_\gamma)$ for all
$\gamma$, and 
\begin{equation}
|\tilde g(t)| \ll \e^{-\sigma |t|} , \qquad \forall t > 0 .
\end{equation}
\end{remark}

\begin{remark}\label{rem2}
We may interpret the sum over conjugacy classes in the spirit
of Propositions \ref{propfirst} and \ref{scatter}:
provided $h$ satisfies {\rm(H1), (H2), (H3*)}, we have
\begin{equation}\label{traceEQ2}
\sum_{j=0}^\infty  h(\rho_j)
= \frac{\Area(\scrM)}{4\pi} 
\int_{-\infty}^\infty h(\rho) \tanh(\pi\rho)\,\rho\, d\rho
+ \frac12
\sum_{\gamma\in H_*}  \int_{\RR} h(\rho)\,n_{\Z_\gamma}(\rho)\, d\rho .
\end{equation}
Alternatively, replace the first term on the right hand side in 
\eqref{traceEQ} by \eqref{ditte}, then
\begin{equation}\label{traceEQ3}
\sum_{j=0}^\infty  h(\rho_j)
= -\frac{\Area(\scrM)}{2\pi} \int_0^\infty \frac{g'(t)}{\sinh(t/2)}\, dt
+ \sum_{\gamma\in H_*}  \sum_{n=1}^\infty
\frac{\ell_{\gamma}\,g(n\ell_{\gamma})}{2\sinh(n\ell_{\gamma}/2)}  .
\end{equation}
\end{remark}

\section{The heat kernel and Weyl's law \label{secWeyl}}

As a first application of the trace formula, we now prove 
Weyl's law for the asymptotic number of eigenvalues $\lambda_j$
below a given $\lambda$,
\begin{equation}
N(\lambda)=\# \{ j: \lambda_j \leq \lambda \} ,
\end{equation}
as $\lambda\to\infty$.

\begin{proposition}[Weyl's law]
\begin{equation}
N(\lambda)\sim\frac{\Area(\scrM)}{4\pi}\, \lambda, \qquad \lambda\to\infty.
\end{equation}
\end{proposition}

\begin{proof}
For any $\beta>0$, the test function 
\begin{equation}
h(\rho)=\e^{-\beta \rho^2}
\end{equation}
is admissible in the trace formula. The Fourier transform is
\begin{equation}
g(t)=\frac{\e^{-t^2/(2\beta)}}{\sqrt{4\pi\beta}} ,
\end{equation}
and so \eqref{traceEQ} reads in this special case 
(with $\lambda_j=\rho_j^2+\tfrac14$)
\begin{multline}\label{traceEQheat}
\sum_{j=0}^\infty  \e^{-\beta\lambda_j}
= \frac{\Area(\scrM)}{4\pi} 
\int_{-\infty}^\infty \e^{-\beta(\rho^2+\tfrac14)} 
\tanh(\pi\rho)\,\rho\, d\rho\\
+\frac{\e^{-\beta/4}}{\sqrt{4\pi\beta}}
 \sum_{\gamma\in H_*}  \sum_{n=1}^\infty
\frac{\ell_{\gamma}\, 
\e^{-(n\ell_\gamma)^2/(2\beta)}}{2\sinh(n\ell_{\gamma}/2)} .
\end{multline}
The sum on the right hand side clearly tends to zero in the limit $\beta\to 0$.
Since $\tanh(\pi\rho)=1+O(\e^{-2\pi|\rho|})$ for all $\rho\in\RR$,
we obtain
\begin{equation}
\sum_{j=0}^\infty  \e^{-\beta\lambda_j}
= \frac{\Area(\scrM)}{4\pi \beta} +O(1), \qquad \beta\to 0. 
\end{equation}
The Proposition now follows from a classical Tauberian theorem
\cite{Widder41}.
\qed
\end{proof}

The sum $\sum_{j=0}^\infty  \e^{-\beta\lambda_j}$ represents 
of course the trace
of the {\em heat kernel} $\e^{\beta\Delta}$.

\begin{corollary}\label{tracecor}
The condition {\rm (H3*)} in Theorem \ref{trace} 
and Remark \ref{rem2} can be replaced by {\rm (H3)}.
\end{corollary}

\begin{proof}
Weyl's law implies that, for any $\delta>0$
\begin{equation}
\sum_{j=0}^\infty (1+\lambda_j)^{-1-\delta/2} < \infty ,
\qquad \text{i.e.,} \qquad 
\sum_{j=0}^\infty (1+\Re \rho_j)^{-2-\delta} < \infty .
\end{equation}
To prove this claim, note that
\begin{align}
\sum_{j=0}^\infty (1+\lambda_j)^{-1-\delta/2}
& =
\int_0^\infty (1+x)^{-1-\delta/2} dN(x) \\
& =
(1+x)^{-1-\delta/2} N(x)\bigg|_{x=0}^\infty \\
& \quad + \left(1+\frac{\delta}{2}\right) 
\int_0^\infty (1+x)^{-2-\delta/2} N(x)\, dx
\end{align}
(use integration by parts) which is finite since $N(x)$ grows linearly
with $x$.

If $h$ satisfies (H1), (H2), (H3), then
the function $h_\epsilon(\rho)=h(\rho)\e^{-\epsilon\rho^2}$
clearly satisfies (H1), (H2), (H3*) for any $\epsilon>0$, 
with the additional uniform bound
\begin{equation}
|h_\epsilon(\rho)|\ll (1+|\Re\rho|)^{-2-\delta}.
\end{equation}
where the implied constant is independent of $\epsilon$.
By repeating the calculation that leads to \eqref{gbound},
we obtain the following estimate for the Fourier transform of $h_\epsilon$,
\begin{equation}\label{gbound2}
|g_\epsilon(t)| \ll \e^{-\sigma|t|} ,
\end{equation}
where the implied constant is again independent of $\epsilon$.
Theorem \ref{trace} yields
\begin{equation}
\sum_{j=0}^\infty  h_\epsilon(\rho_j)
= \frac{\Area(\scrM)}{4\pi} 
\int_{-\infty}^\infty h_\epsilon(\rho) \tanh(\pi\rho)\,\rho\, d\rho
+ \sum_{\gamma\in H_*}  \sum_{n=1}^\infty
\frac{\ell_{\gamma}\,g_\epsilon(n\ell_{\gamma})}{2\sinh(n\ell_{\gamma}/2)} .
\end{equation}
Due to the above $\epsilon$-uniform bounds, 
both sides of the trace formula converge
absolutely, uniformly for all $\epsilon>0$. We may therefore take
the limit $\epsilon\to 0$ inside the sums and integral.
\end{proof}

\section{The density of closed geodesics \label{secHuber}}

In the previous section we have used the trace formula to 
obtain Weyl's law on the distribution of eigenvalues $\lambda_j$.
By using the appropriate test function, one can similarly 
work out the asymptotic number of primitive closed geodesic 
with lengths $\ell_\gamma\leq L$,
\begin{equation}
\Pi(L)=\#\{ \gamma\in H_* : \ell_\gamma \leq L \} .
\end{equation}
In view of Remark \ref{rem} we know that, for any $\delta>0$,
\begin{equation}
\sum_{\gamma\in H_*} \ell_\gamma \e^{-\ell_\gamma(1+\delta)} < \infty
\end{equation}
which implies that, for any $\epsilon>0$,
\begin{equation}\label{rough}
\Pi(L) \ll_\epsilon \e^{L(1+\epsilon)} .
\end{equation}
There is in fact an a priori geometric argument (cf. \cite{HejhalI})
which yields this bounds with $\epsilon=0$, but the rough estimate
\eqref{rough} is sufficient for the following argument.

Let us consider the density of closed geodesics
in the interval $[a+L,b+L]$ where $a$ and $b$ are fixed and $L\to\infty$.
To avoid technicalities, we will here only use smoothed
counting functions 
\begin{equation}
\sum_{\gamma\in H_*} \psi_L(\ell_\gamma)
=
\int_0^\infty  \psi_L(t)\, d\Pi(t)
\end{equation}
where $\psi_L(t)=\psi(t-L)$ and $\psi\in\C_0^\infty(\RR)$.
One may think of $\psi$ as a smoothed characteristic function of $[a,b]$.
Stronger results for true counting functions require
a detailed analysis of Selberg's zeta function, which will be introduced 
in Section \ref{selzeta}.

Let $\rho_0,\ldots,\rho_M$ be those $\rho_j$ with $\Im \rho_j < 0$.
The corresponding eigenvalues $\lambda_0,\ldots,\lambda_M$ are
referred to as the {\em small} eigenvalues.

\begin{proposition}\label{pgt}
Let $\psi\in\C_0^\infty(\RR)$. Then, for $L>1$,
\begin{equation}
\int_0^\infty  \psi_L(t)\, d\Pi(t)
=
\int_0^\infty  \psi_L(t)\, d\widetilde\Pi(t) 
+ O\left(\frac{\e^{L/2}}{L}\right),
\end{equation}
where
\begin{equation}
d\widetilde\Pi(t) = \sum_{j=0}^M \frac{\e^{(\tfrac12+\i\rho_j) t}}{t} \, dt .  
\end{equation}
\end{proposition}

\begin{proof}
The plan is to apply the trace formula with
\begin{equation}
g(t)= \frac{2\sinh(t/2)}{t} \left[\psi(t-L)+\psi(-t-L)\right]
\end{equation}
which is even and, for $L$ large enough, in $\C_0^\infty(\RR)$. 
Hence its Fourier transform,
\begin{align}
h(\rho) & = 
\int_\RR \frac1t 
\big(\e^{(\tfrac12 +\i\rho) t} - \e^{(-\tfrac12 +\i\rho) t}\big)
\left[\psi(t-L)+\psi(-t-L)\right]\, dt  \\
& = 
\int_\RR \frac1t 
\big(\e^{(\tfrac12 +\i\rho) t} - \e^{(-\tfrac12 +\i\rho) t}
+ \e^{(-\tfrac12-\i\rho) t} - \e^{(\tfrac12 - \i\rho) t} \big) \psi(t-L)\, dt ,
\end{align}
satisfies (H1), (H2), (H3). Let us begin with the integral
\begin{equation}
\int_\RR \frac1t \,
\e^{(\tfrac12 +\i\rho) t} \psi(t-L)\, dt 
=
\e^{(\tfrac12 +\i\rho)L}
\int_\RR \frac{1}{t+L} \,
\e^{(\tfrac12 +\i\rho) t} \psi(t)\, dt . 
\end{equation}
Repeated integration by parts yields the upper bound
\begin{align}
\bigg| \int_\RR \frac{1}{t+L} \,
\e^{(\tfrac12 +\i\rho) t} \psi(t)\, dt \bigg|
& \ll_N
\frac{1}{(1+|\rho|)^N}
\int_{\supp\psi} \frac{1}{t+L} \,
\e^{t} \, dt \\
& \ll_N
\frac{1}{L (1+|\rho|)^N} .
\end{align}
So
\begin{equation}
\int_\RR \frac1t \,
\e^{(\tfrac12 +\i\rho) t} \psi(t-L)\, dt 
\ll_N \frac{\e^{(\tfrac12 -\Im\rho)L}}{L (1+|\rho|)^N} .
\end{equation}
This bound is useful for $\Im\rho=0$.
The other corresponding integrals can be estimated in a similar way,
to obtain the bounds (assume $-1/2\leq\Im\rho\leq 0$)
\begin{equation}
\int_\RR \frac1t \,
\e^{(-\tfrac12 +\i\rho) t} \psi(t-L)\, dt 
\ll_N \frac{1}{L (1+|\rho|)^N} ,
\end{equation}
\begin{equation}
\int_\RR \frac1t \,
\e^{(-\tfrac12 -\i\rho) t} \psi(t-L)\, dt 
\ll_N \frac{1}{L (1+|\rho|)^N} ,
\end{equation}
and
\begin{equation}
\int_\RR \frac1t \,
\e^{(\tfrac12 -\i\rho) t} \psi(t-L)\, dt 
\ll_N \frac{\e^{L/2}}{L (1+|\rho|)^N} .
\end{equation}
Therefore, using the above bound with $N=3$, say, yields 
\begin{equation}
\sum_{j=0}^M  h(\rho_j) =
\int_0^\infty  \psi_L(t)\, d\widetilde\Pi(t) + 
O\left(\frac{\e^{L/2}}{L}\right)
\end{equation}
and
\begin{equation}
\sum_{j=M+1}^\infty  h(\rho_j)
- \frac{\Area(\scrM)}{4\pi} 
\int_{-\infty}^\infty h(\rho) \tanh(\pi\rho)\,\rho\, d\rho \\
=
O_N\left(\frac{\e^{L/2}}{L}\right) .
\end{equation}
The sum of the above terms equals, by the trace formula, the expression
\begin{equation}
\sum_{\gamma\in H_*} \sum_{n=1}^\infty \frac1n  \psi_L(n\ell_\gamma) .
\end{equation}
The a priori bound \eqref{rough} tells us that terms with $n\geq 2$
(corresponding to repetitions of primitive closed geodesics)
are of lower order. To be precise,
\begin{equation}
\sum_{\gamma\in H_*} \sum_{n=2}^\infty \frac1n  \psi_L(n\ell_\gamma) 
\ll_\epsilon
\sum_{2\leq n \leq (L+b)/\ell_{\min} } \frac1n  \,  \e^{(L+b)(1+\epsilon)/n}
\ll_\epsilon
\e^{L(1+\epsilon)/2} ,
\end{equation}
where we assume that $\psi_L$ is supported in $[a+L,b+L]$,
and $\ell_{\min}$ is the length of the shortest primitive closed geodesic.
Therefore
\begin{equation}\label{rro}
\sum_{\gamma\in H_*} \psi_L(\ell_\gamma)
=
\int_0^\infty  \psi_L(t)\, d\widetilde\Pi(t) 
+ O\left(\e^{L(1+\epsilon)/2}\right) .
\end{equation}
The leading order term as $L\to\infty$ is
\begin{equation}
\sum_{\gamma\in H_*} \psi_L(\ell_\gamma)
\sim
\int_0^\infty  \psi_L(t)\, \frac{\e^t}{t}\, dt
\ll
\frac{\e^{L+b}}{L+b}
\end{equation}
which leads to the improved upper bound for the sum involving repetitions,
\begin{equation}
\sum_{\gamma\in H_*} \sum_{n=2}^\infty \frac1n  \psi_L(n\ell_\gamma) 
\ll
\sum_{2\leq n \leq (L+b)/\ell_{\min} } \frac{\e^{(L+b)/n}}{L+b}
\ll
\frac{\e^{L/2}}{L} ,
\end{equation}
and hence leads to the desired improved error estimate in \eqref{rro}.
\qed
\end{proof}

\section{Trace of the resolvent\label{secResolvent}}

The trace of the resolvent $R(\lambda)=(\Delta+\lambda)^{-1}$
is formally 
\begin{equation}
\Tr R(\lambda)= \sum_{j=0}^\infty (\lambda-\lambda_j)^{-1}
= \sum_{j=0}^\infty h(\rho_j), 
\qquad h(\rho')= (\rho^2-{\rho'}^2)^{-1} ,
\end{equation}
where $\rho=\sqrt{\lambda-\tfrac14}$ as usual.
The test function $h$ does not, however, respect condition (H3). 
To overcome this difficulty, we define the regularized resolvent
\begin{equation}
\widetilde R(\lambda)=(\Delta+\lambda)^{-1}-(\Delta+\lambda_*)^{-1}
\end{equation} 
for some fixed $\lambda_*$.
The corresponding test function is 
\begin{equation}
h(\rho')= (\rho^2-{\rho'}^2)^{-1}-(\rho_*^2-{\rho'}^2)^{-1},
\end{equation}
which clearly satisfies (H3), since
\begin{equation}
h(\rho')= \frac{\rho_*^2-\rho^2}{(\rho^2-{\rho'}^2)(\rho_*^2-{\rho'}^2)}
= O({\rho'}^{-4}).
\end{equation}
We have already encountered the kernel of the regularized resolvent,
\begin{equation}
k(z,w) = G_\rho(z,w)- G_{\rho_*}(z,w) ,
\end{equation}
in Section \ref{ghost}.
The trace of the regularized resolvent is thus
\begin{equation}
\Tr \widetilde R(\rho) =
\sum_{j=0}^\infty
\left[(\rho^2-{\rho_j}^2)^{-1}-(\rho_*^2-{\rho_j}^2)^{-1}\right] 
\end{equation}
which, for any fixed $\rho^*\notin\{\pm \rho_j\}$, is
a meromorphic function in $\CC$ with simple poles at $\rho=\pm \rho_j$.
$h$ is analytic in the strip $|\Im\rho'|\leq\sigma$
provided $\sigma<|\Im\rho|<|\Im \rho^*|$ where $\sigma>1/2$.

The trace formula \eqref{traceEQ2} implies therefore
(use formula \eqref{ztow} for the first term on the right hand side,
and \eqref{nnn3} for the second)
\begin{multline}
\Tr \widetilde R(\rho)
= - \frac{\Area(\scrM)}{2\pi\i}  
\sum_{l=0}^\infty  
\left[\frac{1}{\rho-\i(l+\tfrac12)} - \frac{1}{\rho_*-\i(l+\tfrac12)}\right]\\
+ \frac{\pi\i}{2\rho} \sum_{\gamma\in H_*} n_{\Z_\gamma}(\rho) + C(\rho_*).
\end{multline}
where
\begin{equation}
C(\rho_*)=
-\frac{\pi\i}{2\rho_*} \sum_{\gamma\in H_*} n_{\Z_\gamma}(\rho_*)
\end{equation}
converges absolutely, cf.~Remark \ref{rem}. 

Let us rewrite this formula as
\begin{multline}\label{rrw}
\frac{1}{2\rho} \sum_{\gamma\in H_*} n_{\Z_\gamma}(\rho)
=
\frac{1}{\pi\i} \sum_{j=0}^\infty
\left[\frac{1}{\rho^2-{\rho_j}^2}
-\frac{1}{\rho_*^2-{\rho_j}^2}\right] \\
-\frac{\Area(\scrM)}{2\pi^2}
\quad \sum_{l=0}^\infty  
\left[\frac{1}{\rho-\i(l+\tfrac12)} 
- \frac{1}{\rho_*-\i(l+\tfrac12)}\right]
-\frac{1}{\pi\i} C(\rho_*).
\end{multline}
All quantities on the right hand side are meromorphic for all $\rho\in\CC$,
for every fixed $\rho_*\in\CC$ away from the singularities
(this is guaranteed for $|\Im\rho_*|>1/2$).
Therefore \eqref{rrw} provides a meromorphic continuation
of
\begin{equation}
n_\Gamma(\rho)
:=\sum_{\gamma\in H_*} n_{\Z_\gamma}(\rho)
\end{equation}
to the whole complex plane.

\begin{proposition} \label{prezeta}
The function $n_\Gamma(\rho)$
has a meromorphic continuation to the whole complex plane,
with
\begin{enumerate}
\item[{\rm (i)}]
simple poles at $\rho=\rho_0,\pm \rho_1,\pm \rho_2,\ldots$ with residue $\res_{\rho_0} n_\Gamma=\frac{1}{\pi\i}$, and 
\begin{equation}
\res_{\pm \rho_j} n_\Gamma = 
\begin{cases}
\frac{2\mu_j}{\pi\i} & \text{ if $\rho_j = 0$,} \\
\frac{\mu_j}{\pi\i} & \text{ if $\rho_j\neq 0$,}
\end{cases}
\end{equation}
where $\mu_j$ is the multiplicity of $\rho_j$.
\item[{\rm (ii)}]
simple poles at $\rho=\i(l+\tfrac12)$ with residue 
\begin{equation}
\res_{\i(l+\tfrac12)} n_\Gamma = 
\begin{cases}
\frac{\Area(\scrM)}{2\pi^2\i}\, (2l+1) + \frac{1}{\pi\i} & \text{ if $l=0$,}\\
\frac{\Area(\scrM)}{2\pi^2\i}\, (2l+1) & \text{ if $l>0$,}
\end{cases}
\end{equation}
\item[{\rm (iii)}]
the functional relation
\begin{equation}\label{funky}
n_\Gamma(\rho)+ n_\Gamma(-\rho) = -\frac{\Area(\scrM)}{\pi}\, 
\rho \tanh(\pi\rho).
\end{equation}
\end{enumerate}
\end{proposition}

\begin{proof}
(i) and (ii) are clear. (iii) follows from the identity \eqref{here}.
\qed
\end{proof}

\section{Selberg's zeta function}\label{selzeta}

Selberg's zeta function is defined by
\begin{equation}
Z(s) = \prod_{\gamma\in H_*} \prod_{m=0}^\infty 
\left( 1-\e^{-\ell_\gamma(s+m)} \right),
\end{equation}
which converges absolutely for $\Re s > 1$; this will become clear below,
cf.~\eqref{logdiv0} and \eqref{logdiv}.
Each factor
\begin{equation}
\prod_{m=0}^\infty \left( 1-\e^{-\ell_\gamma(s+m)} \right)
\end{equation}
converges for all $s\in\CC$, with zeros at
\begin{equation}
s=s_{\nu m}= -m + \i (2\pi/\ell_\gamma) \nu, 
\qquad \nu\in\ZZ, \quad m=0,1,2,\ldots .
\end{equation}
Note that $s_{\nu m}=\tfrac12 + \i \rho_{\nu m}$ where
$\rho_{\nu m}$ are the scattering poles for the
hyperbolic cylinder $\Z_\gamma\backslash\HH^2$.
What is more,
\begin{equation}\label{logdiv0}
\frac{d}{ds} \log \prod_{m=0}^\infty \left( 1-\e^{-\ell_\gamma(s+m)} \right)
=
-\ell_\gamma \sum_{m=0}^\infty \left( 1-\e^{\ell_\gamma(s+m)} \right)^{-1}
= \pi  n_{\Z_\gamma}(\rho), 
\end{equation}
with $s=\tfrac12 + \i \rho$, and thus
\begin{equation} \label{logdiv}
\frac{Z'}{Z}(s) 
= \pi n_\Gamma(\rho),
\qquad \Re s > 1 \quad \text{(i.e., $\Im\rho<-1/2$)}.
\end{equation}

Recall the the genus is related to the area of $\scrM$ by
$\Area(\scrM)=4\pi(g-1)$.

\begin{theorem} \label{zetathm}
The Selberg zeta function can be analytically continued to an entire
function $Z(s)$ whose zeros are characterized as follows.
(We divide the set of zeros into two classes, {\em trivial} and 
{\em non-trivial}).
\begin{enumerate}
\item[{\rm (i)}]
The non-trivial zeros
of $Z(s)$ are located at $s=1$ and 
$s=\tfrac12 \pm\i \rho_j$ $(j=1,2,3,\ldots)$ with multiplicity 
\begin{equation}
\begin{cases}
2\mu_j & \text{ if $\rho_j=0$}\\
\mu_j & \text{ if $\rho_j\neq 0$.}
\end{cases}
\end{equation}
The zero at $s=1$ (corresponding to $j=0$) has multiplicity 1.
\item[{\rm (ii)}]
The trivial zeros are located at at $s=-l$, $l=0,1,2,\ldots$ 
and have multiplicity $2g-1$ for $l=0$ and $2(g-1)(2l+1)$ for $l>0$.
\end{enumerate}
Furthermore $Z(s)$ satisfies the functional equation
\begin{equation}
Z(s)=Z(1-s) \exp \left[ 4\pi (g-1) \int_0^{s-\tfrac12} 
v\tan(\pi v)\,dv\right] .
\end{equation}
\end{theorem}

\begin{proof}
Equation \eqref{rrw} yields
\begin{multline}\label{rrws}
\frac{1}{2s-1}\frac{Z'}{Z}(s)
=
\sum_{j=0}^\infty
\left[\frac{1}{(s-\tfrac12)^2+{\rho_j}^2}
+ \frac{1}{\rho_*^2-{\rho_j}^2}\right] \\
- 2(g-1)
\sum_{l=0}^\infty  
\left[\frac{1}{s+l} 
- \frac{1}{l+\tfrac12+\i\rho_*}\right]
+ C(\rho_*) .
\end{multline}
Note that
\begin{equation}
\frac{2s-1}{(s-\tfrac12)^2+{\rho_j}^2}
= 
\frac{1}{s-(\tfrac12+\i\rho_j)}+\frac{1}{s-(\tfrac12-\i\rho_j)},
\end{equation}
hence the corresponding residue is 1.
Furthermore
\begin{equation}
-2(g-1) \frac{2s-1}{s+l} 
\end{equation}
has residue $2(g-1)(2l+1)$ at $s=-l$.
Statements (i) and (ii) are now evident. The functional relation 
follows from \eqref{funky}, which can be written as
\begin{equation}
\frac{Z'}{Z}(s)+\frac{Z'}{Z}(1-s) = 
4\pi(g-1) (s-\tfrac12)\tan[\pi (s-\tfrac12)] .
\end{equation}
Integrating this yields
\begin{equation}
\log Z(s)- \log Z(1-s) = 
4\pi(g-1) \int_0^{s-1/2} v \tan(\pi v)\, dv + c ,
\end{equation}
that is
\begin{equation}
Z(s)/Z(1-s) = 
\exp\left[ 4\pi(g-1) \int_0^{s-1/2} v \tan(\pi v)\, dv + c \right],
\end{equation}
The constant of integration $c$ is determined by setting $s=1/2$.
Notice that the exponential is independent of the path of integration.
\qed
\end{proof}

One important application of the zeta function is a precise asymptotics
for the number of primitive closed geodesics of length less than $L$,
$L\to\infty$; we have \cite{HejhalI}
\begin{equation}\label{pgt2}
\Pi(L)
=
\int_1^L  d\widetilde\Pi(t) 
+ O\left(\frac{\e^{\frac34 L}}{\sqrt L}\right) .
\end{equation}
The error estimate is worse
than in Proposition \ref{pgt}, since we have replaced
the smooth test functions by a characteristic function.
The asymptotic relation \eqref{pgt2} is often 
referred to as {\em Prime Geodesic
Theorem}, due to its similarity with the {\em Prime Number Theorem}.
The proof of \eqref{pgt2} in fact follows the same strategy
as in the Prime Number Theorem, where the Selberg zeta function
plays the role of Riemann's zeta function.

\section{Suggestions for exercises and further reading\label{secSuggestions}}

\begin{enumerate}

\item
{\bf Poisson summation.}

\begin{enumerate}
\item
The Poisson summation formula \eqref{Poisson} reads in higher dimension $d$
\begin{equation}\label{Poisson30}
\sum_{\vecm\in\ZZ^d} f(\vecm)= \sum_{\vecn\in\ZZ^d} \widehat f(\vecn)
\end{equation}
with the Fourier transform
\begin{equation}
\widehat f(\vectau) = \int_{\RR^d} 
f(\vecrho) \, \e^{2\pi\i \vectau\cdot\vecrho} 
d\rho.
\end{equation}
Prove \eqref{Poisson30} for a suitable class of test functions $f$.

\item
Show that \eqref{Poisson30} can be written in the form
\begin{equation}\label{Poisson31}
\sum_{\vecm\in L^*} f(\vecm)= 
\Vol(L\backslash\RR^d) \sum_{\vecn\in L} \widehat f(\vecn)
\end{equation}
where $L$ is any lattice in $\RR^d$ and $L^*$ its dual lattice.

\item Any flat torus can be represented as the quotient $L\backslash\RR^d$,
where the Riemannian metric is the usual euclidean metric.
Show that the normalized eigenfunctions of the Laplacian are 
\begin{equation}
\varphi_{\vecm}(\vecx)= \Vol(L\backslash\RR^d)^{-1/2} 
\e^{2\pi\i \vecm\cdot\vecx}
\end{equation}
for every $\vecm\in L^*$ and work out the corresponding eigenvalues 
$\lambda_j$.

\item
Use \eqref{Poisson31} to derive a trace formula for
$$
\sum_{j=0}^\infty h(\rho_j)
$$
where $\rho_j=\sqrt{\lambda_j}$ (this formula is the famous
Hardy-Voronoi formula, cf.~\cite{Hejhal76}).
\end{enumerate}

\item
{\bf Semiclassics.}

\begin{enumerate}
\item
Show that for $\rho\to\infty$
\begin{equation} \label{statphase}
G_\rho(z,w) = -\frac{1}{2\pi} \sqrt{\frac{\pi}{2\rho\sinh\tau}}\;
\e^{-\i\rho\tau - \i \pi/4} + O(\rho^{-1}), 
\end{equation}
for all fixed $\tau=d(z,w)>0$.
Hint: divide the integral \eqref{intrep} into the ranges
$[\tau,2\tau)$ and $[2\tau,\infty)$. The second range is easily controlled.
For the first range, use the Taylor expansion for $\cosh t$ at
$t=\tau$ to expand the denominator of the integrand.
Relation \eqref{statphase} can also be obtained from the connection of the
Legendre function with the confluent hypergeometric series $F(a,b,c;z)$,
\begin{equation}
Q_\nu(\cosh\tau)
= \sqrt{\pi} \frac{\Gamma(\nu+1)}{\Gamma(\nu+\tfrac32)}
\frac{\e^{-(\nu+1)\tau}}{(1-\e^{-2\tau})^{1/2}}
F\left(\frac12,\frac12,\nu+\frac32; \frac{1}{1-\e^{2\tau}}\right).
\end{equation}

\item
Show that \eqref{statphase} is consistent with
\cite[eq.~(41)]{Gutzwiller89}. Hint: use the $(s,u)$-coordinates
defined in \eqref{su}.

\item 
Compare the Gutzwiller trace formula \cite{Gutzwiller89}, \cite{Combescure99},
with the Selberg trace formula. (Analogues of 
the {\em ghost of the sphere} (Section \ref{ghost})
for more general systems are discussed in \cite{Berry94}.)

\end{enumerate}

\item
{\bf The Riemann-Weil explicit formula.}

\begin{enumerate}
\item
Compare the Selberg trace formula with the Riemann-Weil explicit formula
\cite[eq.~(6.7)]{Hejhal76}, by identifying Riemann zeros with 
the square-root $\sqrt{\lambda_j-\tfrac14}$ of eigenvalues $\lambda_j$
of the Laplacian, and logs of prime numbers
with lengths of closed geodesics. 
See \cite{Berry99} for more on this.

\item
What is the analogue of the ghost of the sphere? 

\end{enumerate}

\item
{\bf Further reading.}
In this course we have discussed Selberg's trace formula in the simplest
possible set-up, for the spectrum of the Laplacian on a compact
surface. The full theory, which is only outlined in Selberg's original paper
\cite{Selberg56}, is developed in great detail in Hejhal's 
lecture notes \cite{HejhalI}, \cite{HejhalII}, 
where the following generalizations are discussed.

\begin{enumerate}

\item
The discrete subgroup $\Gamma$ may contain elliptic elements,
which leads to conical singularities on the surface, and
reflections (i.e., orientation reversing isometries). Technically
more challenging is the treatment of groups $\Gamma$ which contain
parabolic elements. In this case $\GamHH$ is no longer compact,
and the spectrum has a continuous part, cf.~\cite{HejhalII}.

\item
Suppose the Laplacian acts on vector valued functions $f:\HH^2\to\CC^N$ 
which are not invariant under the action of $T_\gamma$, but satisfy
\begin{equation}\label{commchi}
T_\gamma f = \chi(\gamma) f \quad\forall\gamma\in\Gamma
\end{equation}
for some fixed unitary representation $\chi:\Gamma\to \U(N)$.
The physical interpretation of this set-up, in the case $N=1$, 
is that Aharonov-Bohm flux lines thread the holes of the surface.

\item
The Laplacian may act on automorphic forms of weight $\alpha$, which 
corresponds, in physical terms, to the Hamiltonian for 
a constant magnetic field $B$ perpendicular
to the surface, where the 
strength of $B$ is proportional to $\alpha$. See Comtet's article
\cite{Comtet87} for details.

\end{enumerate}

I also recommend Balazs and Voros'
Physics Reports article \cite{Balazs86} and the books by Buser \cite{Buser92},
Iwaniec \cite{Iwaniec02} and Terras \cite{Terras85}, which give a beautiful 
introduction to the theory. 
Readers interested in hyperbolic three-space
will enjoy the book by Elstrodt, Grunewald and Mennicke \cite{Elstrodt98}.
Gelfand, Graev and Pyatetskii-Shapiro \cite{Gelfand90} take
a representation-theoretic view on Selberg's trace formula.

\end{enumerate}

\begin{acknowledgement}
The author has been supported by an EPSRC Advanced Research Fellowship
and the EC Research Training Network (Mathematical Aspects of Quantum Chaos) 
HPRN-CT-2000-00103.
\end{acknowledgement}


\begin{thebibliography}{99}

\bibitem{Balazs86}
N.L.\,Balazs and A.\,Voros, 
Chaos on the pseudosphere, {\em Phys. Rep.} {\bf 143} (1986) 109-240.

\bibitem{Berry94}
M.V.\,Berry and C.J.\,Howls, 
High orders of the Weyl expansion for quantum billiards: 
resurgence of periodic orbits, and the Stokes phenomenon,
{\em Proc. Roy. Soc. London Ser.} A {\bf 447} (1994) 527-555.

\bibitem{Berry99}
M.V.\,Berry and J.P.\,Keating, 
The Riemann zeros and eigenvalue asymptotics,
{\em  SIAM Rev.} {\bf 41} (1999) 236-266.

\bibitem{Buser92}
P.\,Buser, {\em Geometry and spectra of compact Riemann surfaces}, 
{\em Progr. Math.} {\bf 106}, Birkh\"auser Boston, Inc., Boston, MA, 1992.

\bibitem{Buser}
P.\,Buser, Lectures on hyperbolic geometry,
{\em this volume}.

\bibitem{Cartier90}
P.\,Cartier and A.\,Voros, 
Une nouvelle interpr\'etation de la formule des traces de 
Selberg,
{\em The Grothendieck Festschrift}, Vol. II, 1-67, 
{\em Progr. Math.} {\bf 87}, Birkh\"auser Boston, Boston, MA, 1990. 

\bibitem{Combescure99}
M.\,Combescure, J.\,Ralston and D.\,Robert, 
A proof of the Gutzwiller semiclassical trace formula using 
coherent states decomposition, {\em Comm. Math. Phys.} {\bf 202} (1999)
463-480.

\bibitem{Comtet87}
A. Comtet, On the Landau levels on the hyperbolic plane,  
{\em Ann. Physics}  {\bf 173} (1987) 185-209.

\bibitem{Elstrodt98}
J.\,Elstrodt, F.\,Grunewald and J.\,Mennicke, 
{\em Groups acting on hyperbolic space},
Springer-Verlag, Berlin, 1998.

\bibitem{Gelfand90}
I.M.\,Gelfand, M.I.\,Graev and I.I.\,Pyatetskii-Shapiro, 
{\em Representation theory and automorphic functions},
Academic Press, Inc., Boston, MA, 1990
(Reprint of the 1969 edition).

\bibitem{Gutzwiller89}
M.C.\,Gutzwiller, 
The semi-classical quantization of chaotic Hamiltonian systems. 
{\em Chaos et physique quantique} (Les Houches, 1989), 201-250, 
North-Holland, Amsterdam, 1991. 

\bibitem{Hejhal76}
D.A.\,Hejhal,
The Selberg trace formula and the Riemann zeta function,
{\em Duke Math. J.} {\bf 43} (1976) 441-482.

\bibitem{HejhalI}
D.A.\,Hejhal,
{\em The Selberg trace formula for ${\rm PSL}(2,\RR)$}, Vol. 1. 
Lecture Notes in Mathematics {\bf 548}, 
Springer-Verlag, Berlin-New York, 1976.

\bibitem{HejhalII}
D.A.\,Hejhal,
{\em The Selberg trace formula for ${\rm PSL}(2,\RR)$}, Vol. 2. 
Lecture Notes in Mathematics {\bf 1001}, 
Springer-Verlag, Berlin-New York, 1983.

\bibitem{Iwaniec02}
H.\,Iwaniec, 
{\em Spectral methods of automorphic forms}, 2nd ed.,
{\em Graduate Studies in Mathematics} {\bf 53},
AMS, Providence, RI; 
Rev. Mat. Iberoamericana, Madrid, 2002.

\bibitem{Selberg56}
A.\,Selberg, 
Harmonic analysis and discontinuous groups in weakly symmetric 
Riemannian spaces with applications to Dirichlet series,
{\em J. Indian Math. Soc.} {\bf 20} (1956) 47-87. 

\bibitem{Terras85}
A.\,Terras, 
{\em Harmonic analysis on symmetric spaces and applications} I, 
Springer-Verlag, New York, 1985.

\bibitem{Widder41}
D.V.\,Widder, 
{\em The Laplace Transform},
Princeton Math. Series {\bf 6}, Princeton University Press, Princeton, 1941. 

\end{thebibliography}
\end{document}